\pdfoutput=1
\documentclass{article}
\usepackage{arxiv}

\usepackage[utf8]{inputenc} 
\usepackage[T1]{fontenc}    
\usepackage{url}            
\usepackage{booktabs}       
\usepackage{amsfonts}       
\usepackage{nicefrac}       
\usepackage{microtype}      
\usepackage{lipsum}		
\usepackage{graphicx}
\usepackage{natbib}
\setcitestyle{numbers}
\usepackage{doi}

\usepackage[utf8]{inputenc}
\usepackage{lineno,hyperref}
\usepackage{graphicx}
\usepackage{multirow}
\usepackage[table,xcdraw]{xcolor}
\usepackage{amsmath}
\usepackage{amsfonts}
\usepackage[nolist]{acronym}
\usepackage{soul}
\usepackage{arydshln}
\usepackage{subcaption}
\usepackage{comment}
\usepackage{geometry}
\usepackage[super]{nth}
 \usepackage{pgfplots}
\pgfplotsset{compat=newest}
\usepgfplotslibrary{external}
\usepgfplotslibrary{patchplots}
\usetikzlibrary{shapes,arrows}
\usetikzlibrary{positioning,calc}
\usetikzlibrary{decorations.pathreplacing,decorations.markings,math,arrows.meta}
\usetikzlibrary{spy}
\usetikzlibrary{pgfplots.groupplots}
\usepgfplotslibrary{fillbetween}
\usetikzlibrary{patterns}
\newbox{\bigpicturebox}

\DeclareMathOperator{\E}{\mathbb{E}}
\def\equationautorefname~#1\null{Equation~(#1)\null}

\newcommand{\bfi}{\mbox{\boldmath$w$}}
\newcommand{\bX}{\mbox{\boldmath$X$}}
\newcommand{\bR}{\mbox{\boldmath$R$}}
\newcommand{\bF}{\mbox{\boldmath$F$}}
\newcommand{\bL}{\mbox{\boldmath$L$}}
\newcommand{\bSigma}{\boldsymbol{\Sigma}}
\newcommand{\bZ}{\mbox{\boldmath$Z$}}
\newcommand{\bS}{\mbox{\boldmath$S$}}
\newcommand{\bV}{\mbox{\boldmath$V$}}
\newcommand{\bM}{\mbox{\boldmath$M$}}
\newcommand{\bg}{\mbox{\boldmath$g$}}
\newcommand{\bu}{\mbox{\boldmath$u$}}
\newcommand{\bU}{\mbox{\boldmath$U$}}
\newcommand{\bq}{\mbox{\boldmath$q$}}
\newcommand{\bp}{\mbox{\boldmath$p$}}

\newcommand\figureFontSize{9}

\begin{acronym}
\acro{MOR}{Model Order Reduction}
\acro{FOM}{Full Order Model}
\acro{ROM}{Reduced Order Model}
\acro{ROMs}{Reduced Order Models}
\acro{GAN}{Generative adversarial network}
\acro{MACpROM}{Model Assurance Criterion parametric ROM}
\acro{SHM}{Structural Health Monitoring}
\acro{FE}{Finite Element}
\acro{DOF}{degree of freedom}
\acro{POD}{Proper Orthogonal Decomposition}
\acro{SVD}{Singular Value Decomposition}
\acro{MAC}{Modal Assurance Criterion}
\acro{LHS}{Latin Hypercube Sampling}
\acro{VAE}{Variational AutoEncoder}
\acro{cVAE}{conditional Variational AutoEncoder}
\end{acronym}

\title{VpROM: A novel Variational AutoEncoder-boosted Reduced Order Model for the treatment of parametric dependencies in nonlinear systems}

\author{\hspace{1mm}Thomas Simpson \\
	Dept. of Civil, Environmental, and Geomatic Engr.\\
	ETH Z\"{u}rich\\
	 Z\"{u}rich, Switzerland \\
	\texttt{simpson@ibk.baug.ethz.ch} \\
	\And
	\hspace{1mm}Konstantinos Vlachas \\
	Dept. of Civil, Environmental, and Geomatic Engr.\\
	ETH Z\"{u}rich\\
	Z\"{u}rich, Switzerland \\
	\texttt{vlachas@ibk.baug.ethz.ch} \\
 	\And
	\hspace{1mm}Anthony Garland \\
	Computer Sciences and Mechaincs\\
	Sandia National Laboratories\\
	Albuquerque, USA \\
 	\And
	\hspace{1mm}Nikolaos Dervilis \\
	Dynamics Research Group\\
	University of Sheffield\\
	Sheffield, UK \\
 	\And
	\hspace{1mm}Eleni Chatzi \\
	Dept. of Civil, Environmental, and Geomatic Engr.\\
	ETH Z\"{u}rich\\
	Z\"{u}rich, Switzerland \\
}

\renewcommand{\shorttitle}{VpROM}

\begin{document}
\maketitle

\begin{abstract}
	Reduced Order Models (ROMs) are of considerable importance in many areas of engineering in which computational time presents difficulties. Established approaches employ projection-based reduction such as Proper Orthogonal Decomposition, however, such methods can become inefficient or fail in the case of parameteric or strongly nonlinear models. Such limitations are usually tackled via a library of local reduction bases each of which being valid for a given parameter vector. The success of such methods, however, is strongly reliant upon the method used to relate the parameter vectors to the local bases, this is typically achieved using clustering or interpolation methods. We propose the replacement of these methods with a Variational Autoencoder (VAE) to be used as a generative model which can infer the local basis corresponding to a given parameter vector in a probabilistic manner. The resulting VAE-boosted parametric ROM \emph{VpROM} still retains the physical insights of a projection-based method but also allows for better treatment of problems where model dependencies or excitation traits cause the dynamic behavior to span multiple response regimes. 
Moreover, the probabilistic treatment of the VAE representation allows for uncertainty quantification on the reduction bases which may then be propagated to the ROM response.
The performance of the proposed approach is validated on an open-source simulation benchmark featuring hysteresis and multi-parametric dependencies, and on a large-scale wind turbine tower characterised by nonlinear material behavior and model uncertainty.
\end{abstract}


\section{Introduction}
The use of \ac{ROMs} in structural dynamics simulations forms a main ingredient of research revolving around the use of accelerated surrogates for the purpose of Structural Health Monitoring \cite{EBRAHIMIAN2017,ROSAFALCO2021106604}, digital twinning \cite{SOLMAN2022272,EDINGTON2023109971} and uncertainty quantification \cite{Sudret,LUTHEN2023115875}. Reduced order modeling techniques are often categorised in terms of \emph{purely data-driven} methods or \emph{physics-based} methods. \emph{Purely data-driven} methods employ input and output simulations, or even recorded data, of the system of interest to learn the underlying dynamics \cite{vlachas2018a, guo2019data}. The advent of computing power, leading to increasingly deep machine learning architectures, has rendered such methods extremely capable of recreating even complex dynamics \cite{ZHANG2020, brink2021neural}. However, such methods always remain limited by the breadth and quality of the data used to train them \cite{peherstorfer2016data,garland2022feature}. \emph{Physics-based} methods, on the other hand, allow for the creation of structured \ac{ROM} representations, initiating from the equations of motion and projecting these onto a lower dimensional space upon which they can be solved \cite{Carlberg2009,qian2020lift,najera2023structure}. Such a formulation maintains a stronger physics connotation and is, in this sense, often easier to interpret.


Reduction methods, which rely on the principle of projection, which are capable of addressing nonlinear and/or parametric systems, often exploit \ac{POD} \cite{benner2017model}. These methods involve the execution of evaluations of the \ac{FOM} and the use of the output response of these simulations, the so-called snapshots, for determining an appropriate reduction basis \cite{Carlberg2009a,gobat2022reduced}. An alternative to \ac{POD} methods is that of Proper Generalized Decomposition (PGD) \cite{chinesta2013proper}, whilst PGD is in someways inspired by \ac{POD} methods it has the very key difference of being an a-priori method, not requiring any simulation of the \ac{FOM} in order to construct the \ac{ROM} these methods have also had significant success in various dynamical systems \cite{pgdSoftTissue,CHINESTA2011578}.

A straightforward approach to \ac{POD}-based \ac{MOR} consists in constructing a global POD basis, in which data snapshots from simulations carried out at different points in the parameter/phase space are all stacked together. From this collection, a single (global) projection basis is then extracted. This method is widely used and has proven robust performance for a number of applications \cite{agathos2020parametrized,Christensen,sokratia}, however, in the case of nonlinear systems, the \ac{POD} only provides an optimal approximate \emph{linear} manifold \cite{Kerschen2005}. As such, with moderate to large nonlinearities, the respective \ac{POD}-based reduction can become inefficient and require the retention of several modes (at the cost of reduction), or even fail entirely \cite{zimmermann2017parametric}. Similarly, with parametric models, a global \ac{POD} reduction basis can yield poor performance or become computationally inefficient \cite{amsallem2012nonlinear}.

To this end, alternative strategies can be employed as a remedy, relying on
either enriched low-order subspaces or a pool of local \ac{POD} bases.
The first technique is exemplified in \cite{wu2019modal,mahdiabadi2021non}, where the authors make use of \emph{enriched} reduction bases, in which underlying linear modal or vibration modes-based subspaces are enriched with modal derivatives in order to capture moderate geometric non-linearities.
The second approach, on the other hand, relies on a library of pre-assembled, \emph{local} \ac{ROMs}, which are highly successful in approximating localised phenomena \cite{quarteroni2014reduced, kapteyn2022data}.
In this context, \emph{local} \ac{ROMs} can be defined with respect to time, implying the assembly of projection subspaces, and thus \ac{ROMs}, which only capture the dynamics on a certain time window of the full behavior \cite{amsallem2016pebl}.
Thus, each \ac{ROM} of the corresponding library refers to a different time window of the response, establishing locality.
Alternatively, the local nature of these bases may refer to certain regions or subdomains of the input parameter space \cite{morsy2022reduced}, estimated through uniform \cite{Haasdonk2011,rozza2013reduced} or adaptive error estimators \cite{paul2015adaptive}, which perform a model or basis selection operation between the local \ac{ROM}s during model evaluations.
Such local bases can better deal with more heavily nonlinear and parameter-dependent systems, as they enable an indirect form of clustering of the parameter space (including time if needed), thus providing an accurate subspace approximation for the governing equations of motion at any parametric sample \cite{grimberg2021mesh}.

Although this family of schemes can be considered an established pathway, when deriving an actionable \ac{ROM} that serves across a broad range of operational conditions, efficiency can be compromised and is highly sensitive to the basis selection technique from the assembled library of local \ac{ROMs}.
In this context, state of the art  techniques employ clustering \cite{ghavamian2017pod} or interpolation\cite{amsallem2016pebl} operations performed in the proper manifold to maximise precision. 
On the other hand, recent contributions suggest that machine learning-inspired techniques can increase utility and improve performance of \ac{ROMs}, whilst achieving an automated training process \cite{fresca2022pod,cicci2022deep}. 
Inspired by the latter, in this work, we suggest substituting interpolation- or clustering-based schemes with an ML-based generative model, while retaining the projection-based reduction that guarantees domain-wide accuracy. 
Our approach aims to increase efficiency and robustness of the \ac{ROM} by approximating the generalised mapping between parametric inputs and local projection bases, and the resulting \ac{ROMs}, via the use of generative modeling.



Generative models are a group of statistical models that can serve for generating outputs from observed/simulated systems, under unseen initial conditions, loads, or for properties outside those used in the original training set . This can be accomplished via conditioning on a parametric  vector that reflects the characteristics of the system at hand.
Formally, such a generative model learns the joint distribution $P(\bX,\bp)$ of the observed data $\bX$ and the parameter set $\bp$. As such, a generative model learns the distribution of the data itself hence allowing new samples to be drawn for simulation of new (previously unobserved or not simulated) outputs.
More relevantly for our case, a generative model learns the distribution $P(\bX\vert \bp)$, which is the distribution of the data given a certain parameter vector.
The utility of such a generative model is twofold; via learning the joint distribution of reduction bases s and associated model parameters, we can generate the local basis corresponding to any parameter sample, whilst further capturing the uncertainty on this inference.

One popular branch of generative models is deep generative models, which make use of deep neural networks as powerful and flexible nonlinear approximators that are suitable for modelling complex dependencies.
The two most common examples of modern deep generative models are the \ac{GAN} \cite{Goodfellow2014} and the \ac{VAE} \cite{Kingma2014}.
Both architectures have garnered significant interest in a wide range of fields, ranging from the traditional machine learning subjects of computer vision \cite{Goodfellow2014} and natural language processing \cite{Bowman2015,Pagnoni2018}, to the domains of life sciences for novel molecule development \cite{Sanchez-Lengeling2018,Kardurin2017} and de-noising and analysis of Electron-microscopy images \cite{Prakash2020,Rosenbaum}. Recently, diffusion models  have become the state of the art, beating the performance of GANs and VAEs in typical tasks for generative models such as image and video generation \cite{NEURIPS2021_Diff,ho2022video}.
With regards to structural engineering, significant works utilising deep generative models include the application of \ac{GAN}s to nonlinear model analysis \cite{Tsialiamanis2022} and the use of \ac{VAE}s for wind turbine blade fatigue estimation \cite{Mylonas2021}. 

In this work, we tackle the problem of generating local bases at unseen parameter/input values, making use of a \ac{VAE} as a nonlinear generative model.
The \ac{VAE} model was chosen due to its proven ability to learn highly nonlinear manifolds and efficiently work with high dimensional data.
\ac{VAE}s have the additional advantage of estimating uncertainty on the predicted bases.
The \ac{VAE} model was used to replace clustering or interpolation methods previously used for basis generation \cite{ghavamian2017pod,amsallem2016pebl} with the aim to improve accuracy,
tackle high dimensional dependencies, and allow for uncertainty quantification.
The structure of this paper is organised as follows: Section \ref{ReducedOrderModelling} gives a background on parameteric reduced order modelling and the current state of the art regarding projection methods for treating nonlinear parameteric systems. Section \ref{VpROM} then describes the \ac{VAE} model and how it is used in this work to replace the current state of the art  interpolation or clustering methods in the parametric \ac{ROM}. Section \ref{Results} then demonstrates the use of the \ac{VAE}-boosted parametric \ac{ROM} on 2 example problems. Firstly on a three-dimensional shear frame, modeled with Bouc-Wen hysteretic nonlinearities in its joints, with multi-parametric behavior depending on system properties and excitation characteristics. Secondly, the method is demonstrated on a large scale \ac{FE} model of a wind turbine tower undergoing plastic deformation. In this case, the methodology is combined with a hyper-reduction scheme to demonstrate the full potential for reduction in computational time. In both cases, the \ac{VAE} is shown to outperform the state of the art  methods whilst also allowing for uncertainty in the \ac{ROM} prediction to be captured. Section \ref{Conclusions} concludes the paper by summarising the work and results achieved as well as the limitations of the method and by offering perspectives on future developments.

\section{Parametric Reduced Order Modelling} \label{ReducedOrderModelling}

The context of our work is the physics-based reduction of parameterised dynamical systems to derive an equivalent low-order surrogate of a \ac{FOM}, namely \ac{FE} formulations adopted for nonlinear structural dynamics simulations. 
Such reduced representations provide accelerated system evaluations, which are useful for downstream tasks such as structural health monitoring diagnostic and prognostic, and decision support for operation and maintenance planning of engineered systems. In this context, projection-based reduction has been previously used for delivering response estimates \cite{agathos2020parametrized}, as well as for parameter estimation \cite{tatsis2022hierarchical}, or damage localisation and quantification tasks\cite{agathos2022parametric}.

This section first introduces the nonlinear equations of motion governing the problem at hand. 
Then, a projection-based reduction framework is described, along with the additional components needed for the treatment of parametric dependencies, largely following a similar methodology to the available state of the art  approaches \cite{boncoraglio2021active,amsallem2016pebl}.
The efficiency considerations when propagating the dynamics in the low-order formulation are discussed last.

\subsection{Problem Statement} \label{ProblemStatement}
We assume a general nonlinear dynamical system, characterised by the parameter vector $\bp=[p_1,...,p_{\mathrm{k}}]^{\mathrm{T}} \in \Omega \subset \bR^k$, which captures all system- and excitation-relevant parameters. 
Each realisation of $\bp$ corresponds to a unique configuration of the system at hand. 
Thus, the dynamic behavior of such a system is given by the following set of nonlinear governing equations of motion:
\begin{equation}
\bM(\bp)\ddot{\bu}(t) + \bg\left(\bu(t), \dot{\bu}(t), \bp \right) = \bF (t,\bp),
\label{eq:FOM}
\end{equation}
\noindent
where $\bu(t) \in \bR^n$ represents the response of the system in terms of displacements, $\bM(\bp) \in \bR^{n \times n}$ denotes the mass matrix, and  $\bF(t, \bp) \in \bR^{n}$ the external excitation.
The order of the system is expressed by the variable $n$, termed the full-order dimension, which physically represents the size of the coordinate space and thus the number of degrees of freedom in our system.
This variable indirectly functions as a measure of the computational resources required for the model evaluation at full-order dimensions.
Finally, the nonlinear effects are injected in the restoring force term $\bg\left(\bu(t), \dot{\bu}(t), \bp \right) \in\bR^{n}$.
This term potentially encodes complex nonlinear phenomena of different nature, ranging from material nonlinearity to hysteresis or interface nonlinearities, which, in turn, depend on the parameter vector realisation and the response of the system.

\subsection{Projection-based model order reduction}
In our work, we employ a Galerkin projection-based scheme, as described in \cite{Vlachas2021}, to derive an efficient and accurate reduced order representation for the problem at hand, as described in \autoref{ProblemStatement}.
Several alternative methodologies can be found in \cite{benner2017model}.
We opt for a projection-based approach due to its interpretability and its utility for applications such as higher-level \ac{SHM} systems \cite{tatsis2022hierarchical}. 
Specifically, the derived \ac{ROM} delivers a low-order, yet still physics-based, representation of the full physical space of the model.
Thus, the \ac{ROM} is not limited to capturing displacements, but additionally infers stresses, strains, and accelerations at once, rather than deriving a \ac{ROM} for specific elements or only at a few nodes \cite{lee2020model, Simpson2021}.
This ability of the parametric \ac{ROM} allows for an estimation of the \ac{FOM} response at any given physical field of interest.

The derivation approach of the parametric \ac{ROM} is described in what follows in a step-wise manner. The approach assumes the availability of a high fidelity \ac{FOM}, in our case a \ac{FE} model that spatially discretises the full-order representation of the system in \autoref{eq:FOM}.
Typically, a projection-based \ac{ROM} relies on the premise that the dynamic response, in the present case the solution of \autoref{eq:FOM}, lies in a low-order subspace of size $r$, where $r$ is orders of magnitude smaller than the \ac{FOM} dimension, denoted by $n$ ($r \ll n$).
Thus, the following approximation holds:
\noindent
\begin{equation}
	\bu\left(\bp\right) \approx \bV\left( \bp \right) \bq
\end{equation}
\noindent where $\bV \in \mathbb{R}^{N \times r}$ represents the \ac{ROM} basis that expresses the aforementioned subspace and $\bq \in \mathbb{R}^{r}$ is the respective low-order coordinate vector.
By substituting $\bu$ into \autoref{eq:FOM} and multiplying with $\bu^{T}$, thus performing a Galerkin projection, the following equivalent system is derived:
\begin{equation}
\Tilde{\bM}\left( \bp \right) \ddot\bq\left( t \right) + \Tilde{\bg}\left(\ddot\bq,\dot\bq, t \right) = \Tilde{\bF}\left( \bp, t \right)
\label{eq:ROM}
\end{equation}
where $\Tilde{\bM}=\bV^T\bM\bV$, $\Tilde{\bg}=\bV^T\bg$ and $\Tilde{\bF}=\bV^T \bF$.
Key to a reduction that achieves an accurate low-order representation is the assembly of the projection basis $\bV$.
Following the suggestions in \cite{amsallem2016pebl}, we employ the \ac{POD} technique to this end.
This strategy evaluates \autoref{eq:FOM} for a training set of parameters and harvests response information to form the following matrix:
\noindent
\begin{equation}
	\hat{\bS} = \left[
	\begin{array}{c c c c}
		\hat{\bU}\left( \bp_1 \right) & \hat{\bU}\left( \bp_2 \right) & \ldots & \hat{\bU}\left( \bp_{N_s} \right)
	\end{array}
	\right]
	\label{eq:Snaps}
\end{equation}
\noindent where $\hat{\bS} \in \mathbb{R}^{N \times (N_t \times N_s)}$ is termed as the snapshot matrix, and $\hat{\bU}\left( \bp_i \right) \in \mathbb{R}^{N \times N_t}$ contains the time history of the response for every \ac{DOF} for a given parametric realisation, henceforth termed as a snapshot.
$N_t$ designates the number of simulation time steps, $\bp_i$ is the parametric input for snapshot $i$ and $N_s$ is the number of snapshots.
In turn, the projection basis $\bV$ is assembled via \ac{SVD} of $\hat{\bS}$:
\noindent
\begin{equation}
	\hat{\bS} = \bL \bSigma \bZ^T
    \label{eq:POD}
\end{equation}
\noindent and after truncating $\bL$:
\begin{equation}
	\bV = \left[
	\begin{array}{c c c c}
		\bL_1 & \bL_2 & \ldots & \bL_r
	\end{array}
	\right]
    \label{eq:Vmodes}
\end{equation}
\noindent where $\bL_i$ are columns of matrix $\bL$, termed as modes.
The truncation is applied to obtain the first $r$ principal orthonormal components of the reduction basis $\bV$.
The error measure utilised in \cite{ghavamian2017pod} is herein employed for this purpose.

\subsection{MACpROM: Treatment of parametric dependencies via clustering}\label{MACpROM}

As indicated in the governing equations of motion in \autoref{eq:FOM}, the dynamic behavior of the system depends on a set of parameters that express system properties or traits of the induced excitation. 
Thus, the resulting response of the system is strongly dependent on the parameter vector realisation and may be dominated by localised effects (in the parametric space) due to the corresponding activation of nonlinear terms.
As a result, a large number of truncated modes in \autoref{eq:Vmodes} are required to capture the underlying behavior for the whole parametric domain if the use of a single projection basis is assumed \cite{agathos2020parametrized}. 
Such an approach would lead to a prohibitively large \ac{ROM} dimension, rendering the reduction inefficient or even intractable \cite{zimmermann2017parametric}. 
To this end, an alternative strategy can be employed as a remedy, relying on a pool of local \ac{POD} bases.
Obtaining a pool of \ac{FOM} snapshots and the corresponding \ac{POD} training bases enables the ROM to capture localised effects and utilise proper interpolation or clustering techniques to approximate the response in intermediate points \cite{amsallem2012nonlinear}.

\begin{table}[!ht]
\centering
	\caption{The algorithmic MACpROM framework. The maximum number of clusters can be used instead of a user-defined tolerance during Step 2.}
	\label{tab:MAC}
	\begin{tabular} { p{16cm} }
	    \hline
	    \textbf{Notation:}\\
	     $N_{s}$: Number of training snapshots, $N_{\mathrm{t}}$: Number of simulation timesteps\\
	     $n$: \ac{FOM} dimension, $r$: \ac{ROM} dimension\\
	     \hdashline
	     \textbf{Input:}\\ Parameter vector realisations $\bp=\left[\bp_1, \bp_2, \dots ,\bp_{N_{\mathrm{s}}}\right]$,\\ Realisation $\bp_q$ to be approximated \\
		\textbf{Output:}\\
		Quantity of interest for realisation $\mathbf{p}_{\mathrm{q}}$\\
		\hdashline
		\textbf{Training phase}\\
		$1$:\textbf{for} k=1,...,$N_{s}$ \textbf{do}\\
		\quad $1^a$: Simulate \ac{FOM} for $\bp_{k}$ and obtain $\hat{\bU}\left( \bp_{k} \right) \in \mathbb{R}^{n \times N_{\mathrm{t}}}$ (\autoref{eq:FOM},\autoref{eq:Snaps}). \\
		\quad $1^b$: Obtain local basis $\bV \left( \bp_k \right) \in \mathbb{R}^{n \times r}$ (\autoref{eq:POD},\autoref{eq:Vmodes})\\
		\\
		$2$: Initiate \ac{MAC}-guided sampling and clustering.\\
		\quad $2^a$: Start with an initial sampling rate and an assumption for clustering center(s).\\
		\quad $2^b$: Evaluate \ac{MAC} between local bases of training realisations.\\
		\quad $2^{c}$: Assign each training basis to a cluster, whose basis minimises the respective \ac{MAC}.\\
		\quad $2^{d}$: Identify the maximum obtained MAC and check if exceeds pre-defined tolerance\\ 
		\quad $2^{e}$: If so, define a new cluster center on the point of maximum \ac{MAC}\\
		\quad $2^{f}$: Refine the sampling domain by adding training states between cluster center(s)\\ \quad \quad \quad and maximum \ac{MAC} point(s) and repeat from 1 (if needed).\\
		\\
		\textbf{Evaluation phase}\\
		\quad$1$: Identify cluster $i$ for $\bp_q$ using $k-NN$.\\
		\quad$2$: Assume that $\bV\left( \bp_q \right)$ is approximated using the respective $\bV_{cluster}^i$\\
		\quad$3$: Evaluate the \ac{ROM} of \autoref{eq:ROM} to approximate the quantity of interest\\
		\hline
	\end{tabular}
\end{table}

In our previous work in \cite{vlachas2022coupling}, we have successfully employed a \ac{MAC}-guided clustering scheme, partially following the suggestions in \cite{amsallem2011online} and exploiting a cosine similarity measure, also referred to as a Modal Assurance Criterion in the \ac{SHM} domain.
In this case, the locality on the \ac{POD} bases refers to forming clusters within the original parametric domain. 
Thus, for each parametric sample, the \ac{ROM} utilises a dedicated \ac{POD} basis, termed $\bV_{cluster}^i$, based on the assigned cluster, to accurately reproduce the underlying \ac{FOM} behavior. 

However, since the dynamic behavior in each sample is dominated by localised nonlinear effects, which are captured $\bV_i$, clustering is not performed directly on the parametric samples nor employs the usual distance metrics. 
Instead, an adaptive sampling procedure is followed, exploiting clustering techniques that rely on the truncated modes in \autoref{eq:Vmodes}, resulting from each \ac{FOM} evaluation.
Specifically, the \ac{MAC} is utilised as a comparative measure between projection bases, evaluating their ability to capture similar nonlinear effects.
The \ac{MAC} or vector cosine or cosine similarity, is defined as a scalar constant, expressing an indirect form of confidence when evaluating information originating from different sources \cite{allemang2003modal}.
In our case, the \ac{MAC} serves as a measure of similarity and coherence between truncated modes of neighbouring local bases.
Assuming $\bfi_r$ and $\bfi_s$ correspond to the $i_{th}$ truncated mode from \autoref{eq:Vmodes} for the projection bases $\bV_r$ and $\bV_s$ respectively, 
the mathematical expression for the \ac{MAC} reads:
\begin{equation}
	\mathrm{MAC}({\bfi_r},{\bfi_s})=\frac{|{\bfi_r}^{T}{\bfi_s}|^{2}}{({\bfi_r}^{T}{\bfi_r})({\bfi_s}^{T}{\bfi_s})}
	\label{eq:MAC}
\end{equation}
In turn, this measure is utilised to evaluate the value of new information that each mode captures, thus orienting an adaptive sampling and a subsequent clustering formulation during the training phase.

This approach is termed here as the \ac{MACpROM} for reference purposes and its respective algorithmic framework is summarised in \autoref{tab:MAC}.
The elements of this approach have been validated in previous works and have been shown capable of delivering an accurate and efficient reduced-order representation of nonlinear systems \cite{vlachas2022coupling,vlachas2022parametric}.
In this work, this approach serves as an established reference scenario for the validation of the proposed \ac{VAE}-boosted \ac{ROM}, termed \emph{VpROM}.

\subsection{CpROM: Treatment of parametric dependencies via local Basis Coefficients Interpolation} \label{CpROM}

An alternative formulation for treating parametric dependencies in the context of nonlinear \ac{MOR} has been proposed and verified with respect to state of the art  approaches in \cite{Vlachas2021}.
Specifically, the authors proposed an interpolation approach on the local projection bases relying on the established techniques in \cite{amsallem2016pebl}.
To this end, a two-stage projection was introduced, thus allowing the dependence on parameters $\bp$ to be formulated on a separate level from that of the snapshot procedure or the local subspaces.
Thus, after constructing a pool of local bases, as described previously for the \ac{MACpROM}, each local basis $\bV_i$ is projected to the assembled global \ac{POD} basis of the domain $\bV_{global}$ through a coefficient matrix $ \bX $ as follows:
\begin{equation}
    \bV_i = \bV_{global} \bX_i
    \label{eq:coeffs}
\end{equation}
\noindent
where $\bV_i \in \mathbb{R}^{n \times r}$, $\bV_{global} \in \mathbb{R}^{n \times \tilde{r}}$ and $ \bX \in \mathbb{R}^{\Tilde{r} \times n}$. 
The variable $\Tilde{r}$ signifies the number of modes retained on the global basis $\bV_{global}$.
In this manner, interpolation can be performed on the level of coefficient matrices $\bX_i$, which comprise a reduced dimension ($\Tilde{r} \ll n$), thus removing any dependency on the large \ac{FOM} dimension $n$ and rendering the required operations more efficient.
In turn, the respective matrices $\bX_i$ are projected and interpolated in an element-wise manner on the tangent space of the proper Grassmannian manifold and projected back on the original space to obtain the respective local \ac{ROM} basis $\bV$ for any validation sample. 
This strategy is required for the local bases to retain certain orthogonality and positive-definiteness properties and is described in detail in \cite{amsallem2012nonlinear,zimmermann2017parametric}.
A schematic visualisation of the approach, along with its algorithmic framework can be found in \cite{Vlachas2021}.
In our work, this framework is termed CpROM, adopting the same acronym as in the original work for reference purposes.
This serves as an additional comparison \ac{ROM} framework to validate the performance of the proposed \emph{VpROM}.

\subsection{Hyper-reduction}\label{hyperreduction}
Both the \ac{MACpROM} and the CpROM frameworks, which were previously presented, are employed in conjunction with an additional operation, known as hyper-reduction, in order to achieve a substantial reduction in computational cost when dealing with nonlinear systems. 
Hyper-reduction refers to a second-tier approximation strategy, that addresses the bottleneck of updating and reconstructing the \ac{ROM} system matrices due to the presence of the nonlinear terms in an online manner  \cite{peherstorfer2014localized}.
In essence, this technique relies on a weighted evaluation of the corresponding projections of the nonlinear terms only at a subset of the total elements in the spatial discretisation, thus providing substantially accelerated model evaluations.
The detailed description of the method and the discussion of the existing alternative approaches are already covered in previous works and, thus, lie beyond the scope of this paper. 
The interested reader can refer to \cite{Farhat2015,peherstorfer2020model, peherstorfer2015online}.
The validation case studies in our work make use of the Energy Conserving Mesh Sampling and Weighting (ECSW) technique presented in \cite{farhat2014dimensional,grimberg2021mesh}. 
\section{VpROM: Coupling of Generative Models with projection-
based ROMs} \label{VpROM}
Current state of the art  methods for the creation of low-order surrogates of parameterised nonlinear dynamical systems rely on the use of interpolation or clustering methods for the estimation of local bases for given parametric configurations. This work introduces a nonlinear generative model, exploiting a \ac{cVAE} formulation, in place of these methods, with the aim of improving robustness and performance issues of the \ac{ROM} by allowing for nonlinearities in the parameter-basis relation to be captured and for high dimensional dependencies to be better-dealt with. Furthermore, the derived \emph{VpROM} allows for increased utility with regard to the capturing of uncertainty in the predicted bases.

\subsection{Variational Autoencoder (VAE)}
The \ac{VAE} first described in \cite{Kingma2014} is a latent variable model, that is, a model in which it is assumed that the observations are driven by certain unobserved \emph{latent} variables.
Such latent variable models are popular in many areas of science and engineering and are often used to reduce the effective dimensionality of data since the dimension of the latent variables is reduced compared to that of the observations \cite{Bishop1998}.
Indeed, such a concept is inherent to structural dynamics, as modal analysis, and its nonlinear extensions, all utilise lower dimensional representations to simplify the required analysis.
In a probabilistic sense, modes can be considered latent variables, which are unobserved and drive the observed dynamics of the system. The \ac{VAE} architecture thus serves to infer relationships between the latent variables and the observed variables by means of deep neural network functions.

In the context of \ac{MOR}, a \ac{VAE} can be considered as a Bayesian implementation of a deterministic autoencoder; a popular deterministic deep learning technique, which has been often exploited for dimensionality reduction. 
Lee and Carlberg \cite{LEE2020} make use of a convolutional autoencoder, in conjunction with the nonlinear Galerkin method to construct \ac{ROM}s for advection-dominated problems dynamics problems.
Further work has combined an autoencoder with statistical regression methods to create fully data-based \ac{ROMs} of nonlinear dynamical systems for structural dynamics \cite{Simpson2021, Simpson2022}. A similar methodology known as Learning Effective Dynamics has also been shown to be effective in creating \ac{ROM}s of some more classical nonlinear dynamical systems in various scientific disciplines \cite{vlachas2022multiscale,VlachasP2022}.

In a \ac{VAE} model it is assumed that the data $\bX$ are characterised by a probabilistic distribution $p(\bX)$, which we would like to approximate by means of a parameterised, and possibly simplified, distribution $p_{\phi}(\bX)$ which is parameterised by the vector $\phi$.
We assume that the complex distribution of the data is driven by a lower dimensional and more simply distributed hidden variable set $\bZ$, with assumed prior distribution $p(\bZ)$.
The concept here is that given a sufficiently powerful and flexible approximator, it is possible to learn a function that maps the simply distributed latent variables $\bZ$ to the complexly distributed data $\bX$ by learning the distribution $p_{\phi}(\bX\vert \bZ)$ \cite{Doersch2016}. This approximated distribution is found in the form of a deep neural network, namely the \emph{decoder network}. This network is then parameterised by $\phi$ corresponding to the weights and biases. This results in the following expression for the generative model:

\begin{equation}
    p_{\phi}(\bX) = \int{p_{\phi}(\bX\vert \bZ)p(\bZ)dz}
\end{equation}
The training of such a generative model necessitates the inference of those \emph{decoder} parameters, $\phi$, that maximise the likelihood of the observations. It is here noted that the term observations in this case refers to synthetically generated data from \ac{FOM} snapshots.
This can be expressed as:
\begin{equation}
    \phi = \operatorname*{argmax}_\phi \prod_{i=1}^{N}\int{p_{\phi}(\bX^i\vert \bZ)p(\bZ)dz}
\end{equation}

\begin{figure}[h]
    \centering
    \includegraphics[width=0.75\textwidth]{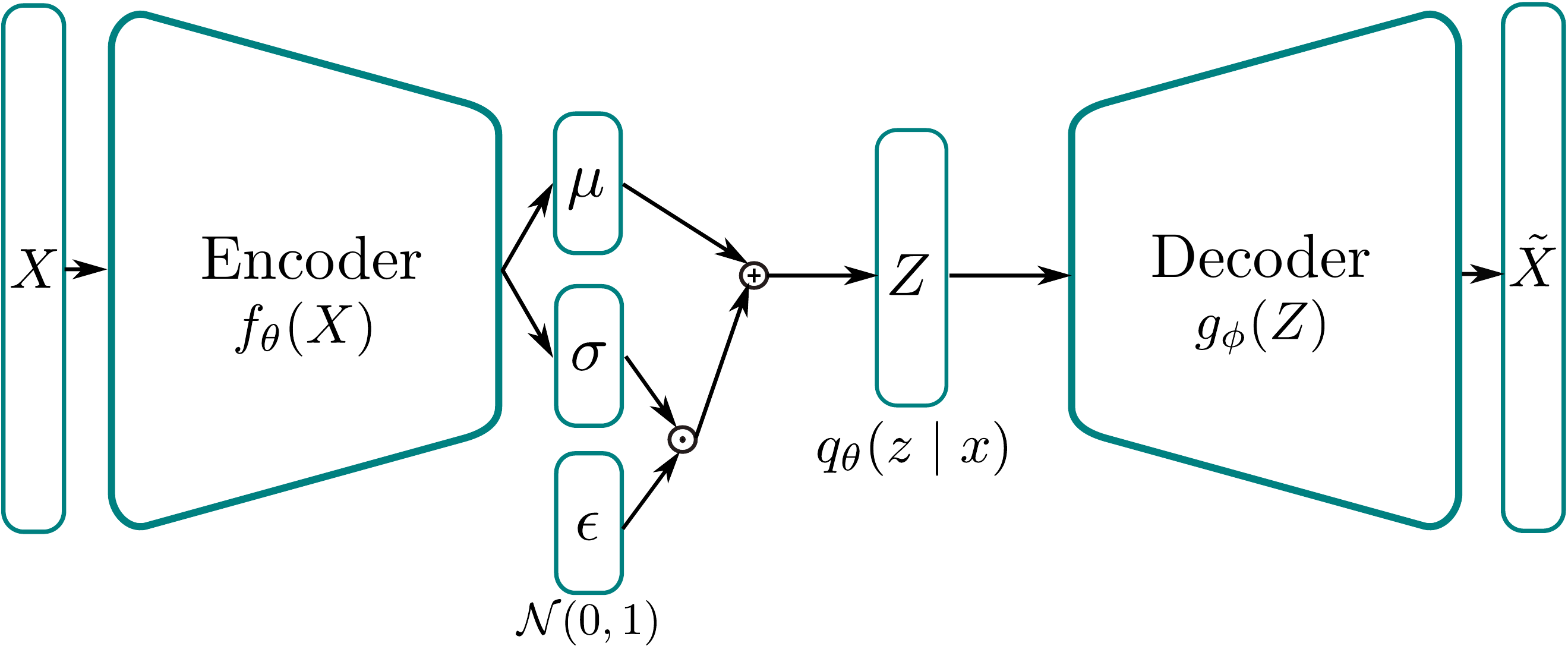}
    \caption{Architecture of a variational autoencoder (VAE).}

    \label{fig:vae}
\end{figure}

The evaluation of this integral, however, presents a problem, as it is generally analytically intractable and computationally inefficient to approximate via sampling.
For this reason, a second \emph{encoder} network is introduced in the typical \ac{VAE} setup, which is additionally parameterised by $\theta$.
This allows the intractable posterior $p(\bZ\vert \bX)$ to be approximated by the parameterised distribution $q_{\theta}(\bZ\vert \bX)$ and hence creates a mapping from the observation space to the latent space. A variational approximation of the distribution of the latent space variable is also made whereby it is assumed that the latent variable takes on a certain known distribution $p(\bZ)$.
This variational approximation of the true posterior results in the following lower bound on the log-likelihood:
\begin{equation}\label{eq:elbo}
    \mathcal{L}(\theta,\phi,\bX)=\E_{q_{\theta}(\bZ\vert \bX)}[log(p_{\phi}(\bX \vert \bZ)]-D_{KL}(q_{\theta}(\bZ\vert \bX)\vert\vert p(\bZ)))
\end{equation}
\noindent
where $D_{KL}$ denotes the Kullback-Leibler divergence, a metric used to measure the similarity of distributions.

It is then this lower bound, known as the evidence based lower bound (ELBO), which is optimised with respect to the parameters of the two networks, $\theta$ and $\phi$. The maximisation of this function aims to i) improve the expected reconstruction loss of the decoder, i.e., the success in recovering observations from the latent variables, and ii) to minimise the KL divergence (and thus maximize the similarity) between the true and approximate posterior of the latent space. Once the VAE is trained based on this process, it is then possible to sample the latent space, using the inferred variational distribution, $q_{\theta}(\bZ\vert \bX)$, and subsequently employ the decoder in order to recreate desired quantities of interest (outputs).

To optimise \autoref{eq:elbo}, it is necessary to estimate the gradients of the ELBO. Kingma et al. \cite{Kingma2014} achieved this by using the \emph{re-parameterisation trick}. By choosing the form of the approximate posterior to be a diagonal Gaussian parameterised by the encoder network, sampling from this distribution can be re-parameterised as follows;

\begin{align}
\eta &= \mathcal{N}(0,I)\\
q_{\theta}(\bZ\vert \bX) &= \mathcal{N}(\bZ:\mu_{\theta}(\bX),\sigma_{\theta}(\bX))\\
& =\mu_{\theta}(\bX)+\eta\odot\sigma_{\theta}(\bX)
\end{align}

 in which $\eta$ represents a sample from the diagonal Gaussian distribution $\mathcal{N}(0,I)$ and $\mu_{\theta}(\bX),\sigma_{\theta}(\bX)$ are the mean and standard deviation values of the latent space as output by the encoder network. The approximate posterior is then reformulated as being a stochastic draw $\eta$ and the deterministic mean and standard deviation values are predicted by the encoder.

This allows for evaluation of the expectation through sampling from a standard multivariate Gaussian, whilst the gradient can be assessed deterministically for each sample allowing the use of back propagation for training. Further, with the choice of a spherical unit Gaussian prior for $p(\bZ)$, the KL divergence term can be calculated analytically \cite{Kingma2014}. This results in the following per sample, differentiable cost function, in which the expectation is evaluated using $N_v$ samples from the latent space per data point.

\begin{equation}
    \mathcal{L}(\theta,\phi,x^i)=\frac{1}{2}\sum^J_{j=1}(1+log((\sigma_j^{(i)})^2)-(\mu_j^{(i)})^2-
    (\sigma_j^{(i)})^2)+\frac{1}{L}\sum_{l=1}^{N_v}log(p_{\phi}(X^i,Z^{i,l})
\end{equation}
where $J$ denotes the dimension of the latent space and $Z^{i,j} =\mu_{\theta}(X^i)+\eta^l\odot\sigma_{\theta}(X^i) $. The number of samples to take from the latent space to evaluate the expectation can even be 1 as in the original formulation of Kingma et al. \cite{Kingma2014}.

As mentioned above, in the VAE model the encoder and decoder functions are approximated using deep neural networks (DNNs). DNNs are a very widely used class of models that exploit multiple neural network layers applied one after another in order to approximate very complex functions more efficiently than shallow networks \cite{BengioLecun}. A thorough description of DNNs and their training can be found in \cite{Goodfellow2014}, in the work herein the DNNs utilised only made use of fully connected layers. In fully connected layers, the transform performed by each layer consists of the matrix multiplication of the input vector with a trainable weights matrix and the addition of a trainable bias vector. A nonlinear activation function, often a $tanh$ or $sigmoid$ is then applied element-wise to the output of this operation. This results in a very flexible and powerful model for learning general nonlinear relations \cite{Goodfellow2014}.

\subsection{VpROM: A conditional Variational Autoencoder (cVAE)-boosted ROM}
\begin{figure}[h]
    \centering
    \includegraphics[width=0.75\textwidth]{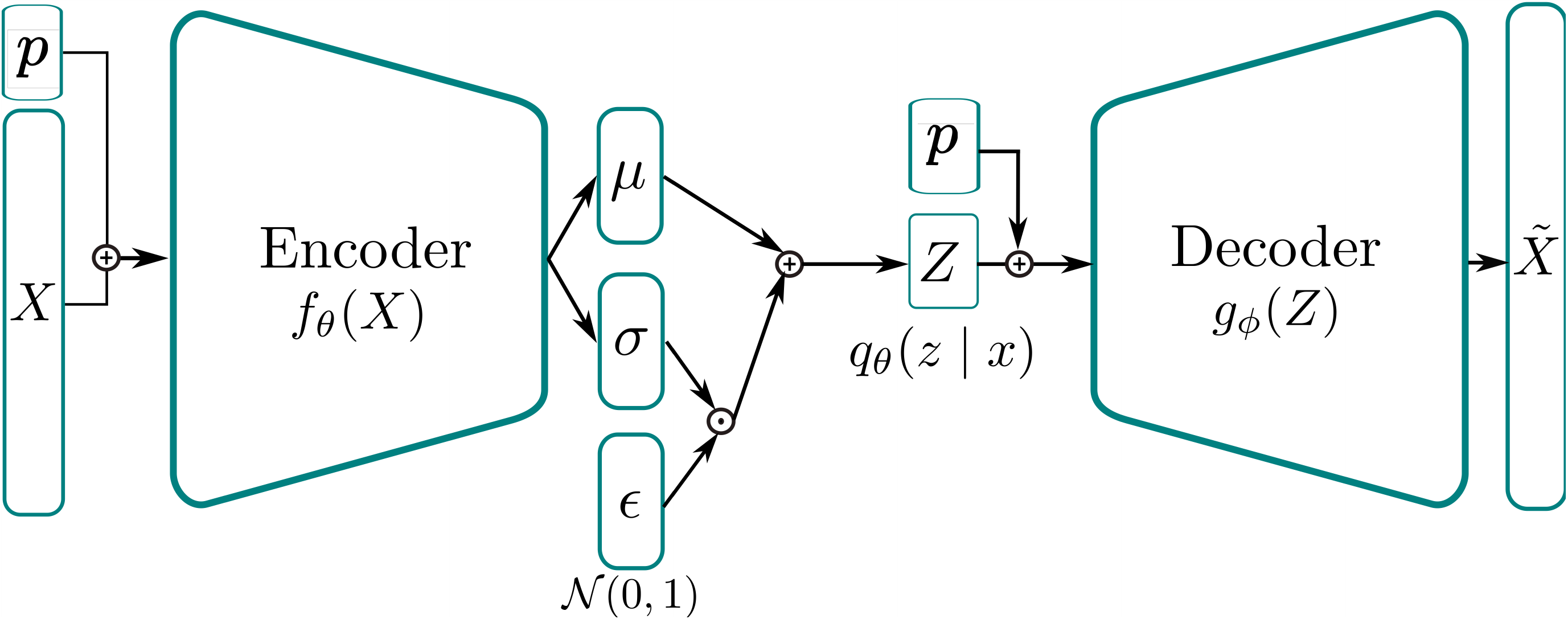}
    \caption{Architecture of a cVAE in which the conditioning variables are injected via concatenation with both the input vector $\bX$ and the latent vector $\bZ$. The input refers to the \ac{ROM} basis coefficients $\bX$ in \autoref{eq:coeffs}}
    \label{fig:cvae}
\end{figure}

In our use case, we don't simply wish to sample possible, plausible bases for our system of interest, rather, we want to sample these bases \emph{conditioned} on given system and load (excitation) parameters. We can achieve this relatively straightforwardly by concatenating the conditioning parameters $\bp$ with the inputs of the \ac{VAE} $\bX$, and with the latent space variables $\bZ$ as demonstrated in \autoref{fig:cvae}. Mathematically, the distribution approximated by the encoder now becomes $q_{\theta}(\bZ \vert \bX,\bp)$ and the distribution approximated by the decoder becomes $p_{\phi}(\bX \vert \bZ,\bp)$. To clarify the role of \ac{cVAE} in the derived \emph{VpROM}, the input referred to in \autoref{fig:cvae} corresponds to the \ac{ROM} basis coefficients $\bX$ in \autoref{eq:coeffs}.
Thus, \autoref{fig:cvae} serves as a visualisation of the mapping process that the \ac{cVAE} carries out to relate the parametric dependencies of the \ac{FOM} with the \ac{ROM} projection basis $\bV$. 
The model dependencies are expressed in \autoref{eq:FOM} and captured in variable $\bp$, whereas the relation to $\bV$ from \autoref{eq:ROM} is expressed through the coefficients $\bX$ from \autoref{eq:coeffs}.

\subsection{VpROM: Generating New Bases}
\begin{figure}[h]
    \centering
    \includegraphics[width=0.375\textwidth]{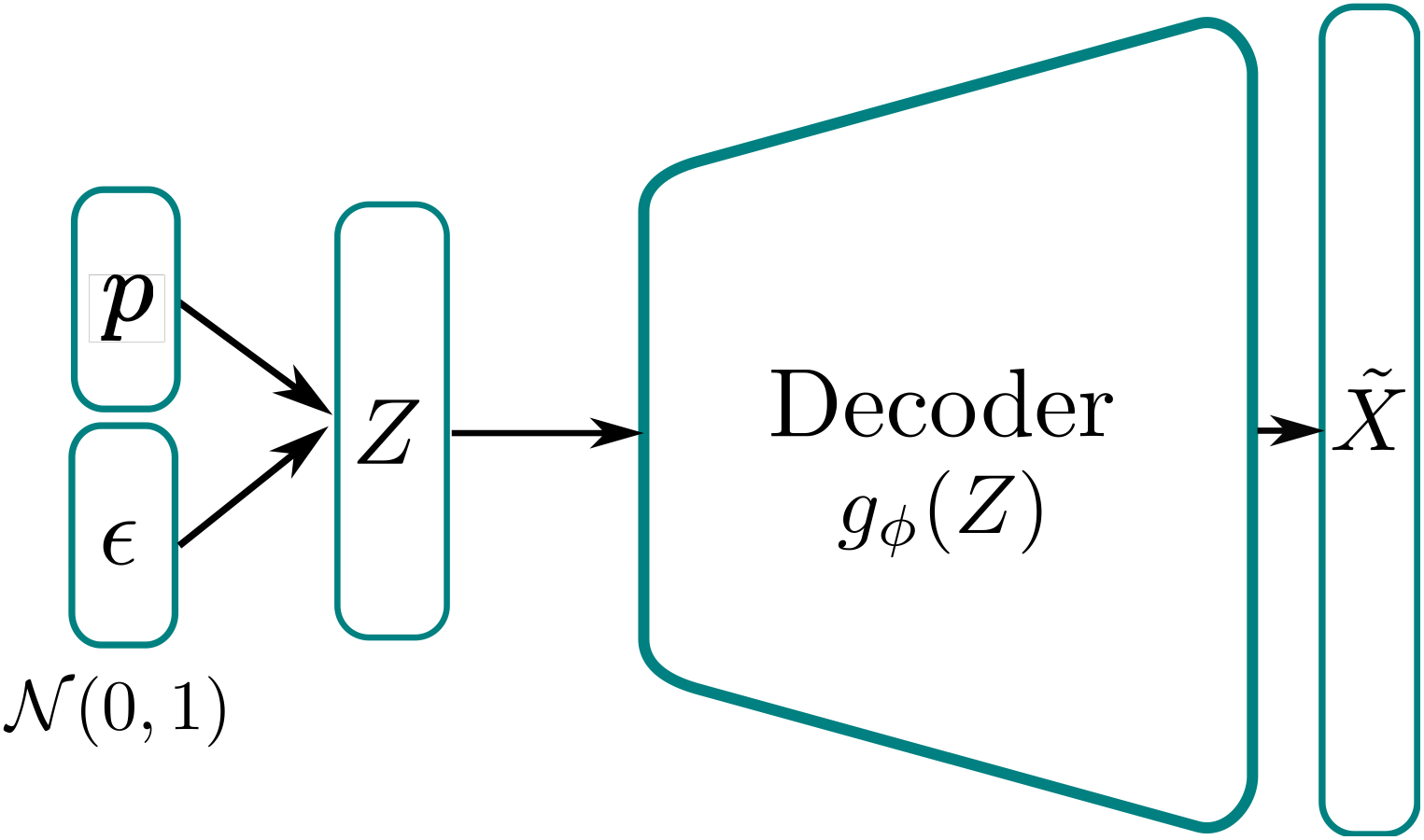}
    \caption{Architecture of the cVAE when used for basis generation: sample latent vectors are taken from the prior distribution $\epsilon$ and concatenated with the conditioning vector $\bp$ before these latent vectors are decoded to find the generated basis coefficients $\bX$ from \autoref{eq:coeffs}.}
    \label{fig:VAE_sample}
\end{figure}

In this work, we propose to train the \ac{cVAE} for creating a generative model, which can be sampled in order to produce the reduced basis coefficients $\bX_i$ from \autoref{eq:coeffs} for a given parameter vector $\bp_i$. In which the parameters either reflect certain properties of the system or of the loading applied. Concretely, once a trained \ac{VAE} is available, the encoder portion is no longer used and predictions are made purely using the decoder and the assumed variational distribution on the latent space. In this case, a diagonal Gaussian is used as shown in \autoref{fig:VAE_sample}. 

Hence, in order to sample from the decoder we simply take a draw, $\epsilon$, from the chosen prior distribution $p(\bZ)$ and concatenate this draw with the given parameters $\bp$. We then pass this latent vector through the decoder, which results in a sample from the observation distribution $p_{\phi}(\bX \vert \bZ,\bp)$. This sample being a single draw from the distribution of the predicted coefficient values for a given parameter vector. Multiple such samples can then be taken by repeating this process in order to find quantities such as the mean and standard deviation of the predicted coefficients. This generative procedure emphasises the importance of minimising the KL divergence term in the loss function. If the KL loss is low then the the approximate posterior distribution $q_{\theta}(\bZ \vert X,p)$ better approaches the prior $p(\bZ)$.

\subsection{Training the VpROM}

As mentioned previously, we wish to train the \ac{cVAE} to generate not the local bases themselves, but rather the coefficient matrices $\bX_i$ as introduced in \autoref{eq:coeffs}, which are then used to generate the local \ac{ROM} subspaces from a global basis. To do this, we require training pairs of parameter vectors $\bp_i$ and corresponding coefficient matrices $\bX_i$. These training pairs are created by first sampling the parameter vector with \ac{LHS}. Each of these training vectors is then used as model/input parameters for a \ac{FOM} snapshot. The generated snapshots are then used to assemble the local bases $\bV_i$ and global basis $\bV_{global}$ respectively, and hence the coefficient matrices, using the procedure described in \autoref{CpROM}. After initial efforts, it was decided that a separate \ac{cVAE} would be trained for each column of the coefficient matrix $\bX_i$. This was considered to be reasonable as each column of $\bX_i$ relates to a different retained \ac{POD} mode, which ought to be mutually orthogonal. Further, the separate consideration of each column results in improved performance of the \ac{cVAE} in generating new bases.

To prepare the data for training, the individual columns $\Tilde{\bX}_{c}$, $c=1,...,n$ of the coefficient matrix $\bX \in \mathrm{R}^{\Tilde{r}*n}$ are taken as vectors and are paired with the parameter vectors $\bp_i$ for training. Further, the parameter values are normalised between $-1$ and $1$ and the coefficient vectors are normalised, as shown in \autoref{eq:coeffNorm}. The coefficient vectors $\Tilde{\bX}_c$ are normalised using the application of a natural logarithm; this offers the advantage of rendering the amplitude of the vector components more similar and hence preventing extreme values from dominating the cost function. The addition of the constant $2$ was required such that all the values were greater than zero in order to avoid an error in the logarithm operation. 

\begin{equation}\label{eq:coeffNorm}
\Tilde{\bX}_c^{norm} =  \ln{(\Tilde{\bX}{_c}+2)}   
\end{equation}

The models were all built and trained using Tensorflow and the Adam algorithm \cite{Kingma2015adam}. In training the \ac{cVAE}, the architecture of the network, the number, and size of layers, and the activation functions used, must also be chosen, these can, however, be treated as hyperparameters and can be optimised according to common methods such as grid search. The network architecture has a massive effect on the expressive power and generalisation of the model. 
All trained models made use of only dense and dropout layers 
Dense layers are a traditional fully connected feedforward neural network layer. Between each of the dense layers, a dropout layer was inserted. Dropout is a technique used for regularising deep neural networks, according to which - during each training update - a certain percentage (the dropout percentage) of activation values of a given layer are set to zero. This technique has been shown to improve generalisation performance of deep neural networks \cite{Srivastava}. As such, the architectural hyperparameters to be optimised for each model include the number of dense layers in the encoder and decoder, the number of neurons in these layers, the activation function used by these layers, and the amount of dropout included between each layer. The size of the bottleneck layer, or in other words, the number of latent variables driving the process is also key for reduction. Noteworthy is, however, that in this case it was found that the application of dropout was not beneficial to the performance; as such, this was not used in any of the finally implemented models.

\section{Results}
\label{Results}
In this section, all aspects relevant to the performance of the proposed framework are validated on two case studies featuring parametric dependencies in system properties and excitation traits. 
The proposed \ac{cVAE}-boosted \ac{ROM} is firstly validated on a nonlinear benchmark simulator of a two-story shear frame featuring hysteretic joints \cite{vlachas2021two}, and then on a larger scale simulation, featuring computational plasticity, which is based on the NREL reference 5-MW wind turbine tower \cite{Jonkman2009}.

As already mentioned in \autoref{ReducedOrderModelling}, we offer a comparison across alternate parametric \ac{ROM} configurations in order to offer a comprehensive discussion on the potential and performance limits of the suggested framework.
The first parametric \ac{ROM} configuration refers to the \ac{MACpROM}, as presented in \autoref{MACpROM}.
This employs a \ac{MAC}-guided clustering approach on the local \ac{POD} bases.
Next, the CpROM presented in \autoref{CpROM} is evaluated, following the local basis coefficients interpolation approach in \cite{Vlachas2021}.
These two \ac{ROMs} are assembled employing existing state-of-the-art approaches and serve comparison purposes. 
The last two parametric ROMs are derived based on the \ac{cVAE} framework proposed here to inject parametric variability in the local projection bases.
We evaluate the performance of the proposed \ac{cVAE}-boosted \ac{ROM}, termed \emph{VpROM}, both with and without the inclusion of hyper-reduction, which is described in \autoref{hyperreduction}.
The notation and configuration of these five schemes are summarised in Table \autoref{tab:pROMS} for reference.  
In what follows we present different case studies comparing the performance of these different schemes.

\begin{table}[!ht]
	\caption{Reference table for compared ROMs.}
	\label{tab:pROMS}
	\centering{
		\begin{tabular} { p{4cm} p{10cm}}
			\hline
			Model Reference Name & Description\\
			\hline
			\ac{FOM} & The full order finite element model  \\
			MACpROM & Pool of local bases and MAC-guided clustering. The approach is described in \autoref{MACpROM}. \\
			CpROM & Formulation of coefficients $\bX$ for each local basis and interpolation, as presented in \autoref{CpROM}.\\
			\emph{VpROM} & The \ac{cVAE}-boosted \ac{ROM} as presented in \autoref{VpROM}\\
            HP-* & The respective * ROM additionally equipped with hyper-reduction\\
			\hline
	\end{tabular}}
\end{table}

Regarding computational resources and timing, the validation simulations of the presented examples are implemented using an in-house built \ac{FE} code, based on the suggestions by \cite{bathe2006finite} and tested on a workstation equipped with an \nth{11} Gen Intel(R) Core(TM) i7-1165G7 processor, running at 2.80GHz, and 32GB of memory.
In addition, the reported computational time is averaged over all \ac{FOM} or \ac{ROM} evaluations of each respective set of configurations (training or testing).
The performance of the various frameworks in terms of reproducing the time history responses of the respective dynamic validation case studies is reported as follows: 

\begin{equation}
    err_q = \dfrac{\sqrt{\sum\limits_{i \in \Tilde{N}_{\text{DOF}}} \sum\limits_{j \in \Tilde{N}_{t}} \left(q_{i}^{j}-\Tilde{q}_{i}^{j}\right)^2}}{\sqrt{\sum\limits_{i \in \Tilde{N}_{\text{DOF}}} \sum\limits_{j \in \Tilde{N}_{t}} \left(q_{i}^{j}\right)^2}} \times 100 \%
    \label{eq:errors}
\end{equation}
\noindent where $\Tilde{N}_{\text{DOF}}$ represents a set of \ac{DOF}s selected for response comparison, $\Tilde{N}_{t}$ a set of selected time steps, $q_{i}$ is the \ac{FOM} quantity of interest at \ac{DOF} $i$, and $\Tilde{q}_{i}$ is the respective inferred value computed using the \ac{ROM} approximation.

\subsection{Two story shear frame with hysteretic links}\label{example:BWlinks}

As an initial example, we consider a \ac{FE} model of a three-dimensional two-story shear frame with nonlinear nodal couplings, each exhibiting a Bouc-Wen hysteretic nonlinearity \cite{ismail2009hysteresis}.
This example is chosen as a demonstrative case study due to the inherent ability of the simulator to model multiple simultaneously activated instances of nonlinearity, thus challenging the precision of any derived \ac{ROM}.
Because this is a low dimensional example, we use it with the main purpose of assessing the accuracy of the respective \ac{ROM}s from \autoref{tab:pROMS}, as the model is rather trivial to demonstrate any substantial computational savings.
The capacity of the various parametric ROMs in terms of accelerating model evaluations is documented in the next case study, featuring a large-scale wind turbine tower.

\begin{figure*}[!ht]
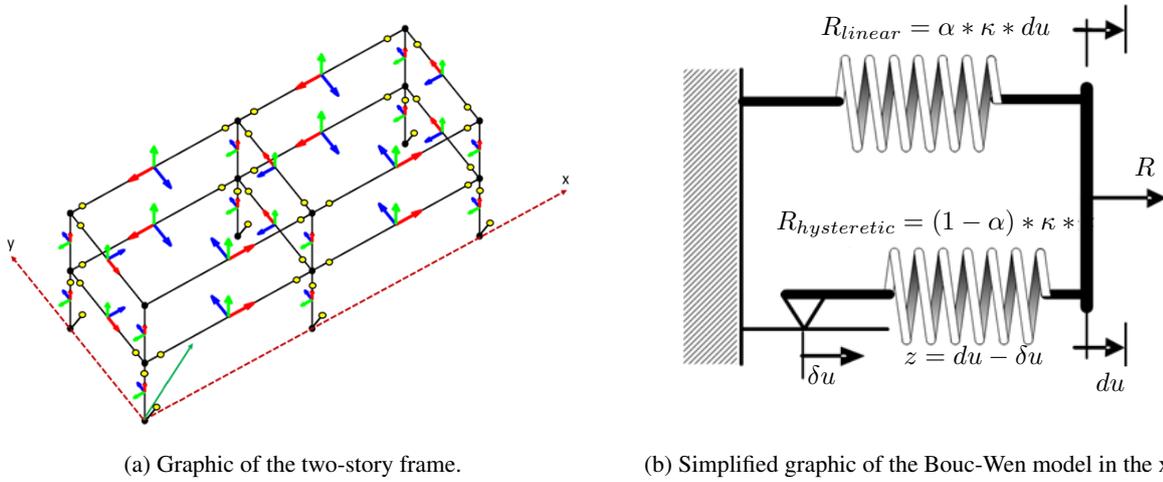

    \centering
    \begin{subfigure}[b]{0.485\textwidth} 
    \input{figures_tikz/Figure4a.tikz}
    \caption{Graphic of the two-story frame. \label{fig:links}}
    \end{subfigure}
    \quad
    \begin{subfigure}[b]{0.485\textwidth}
     \centering
     \input{figures_tikz/Figure4b.tikz}
     \caption{Simplified graphic of the Bouc-Wen model in the x-\ac{DOF}.\label{fig:bw}}
     \end{subfigure}
\caption{Graphic of the frame setup and illustration of the nonlinear mechanism in its links. The green arrow indicates the direction of ground motion, whereas the colored arrows the orientation of the beam elements.  \label{fig:BWmodel}}   
\end{figure*}

A graphical illustration of the setup of the shear frame is visualized in \autoref{fig:links}, where the hysteretic links assume no length, although the virtual nodes are depicted within a distance from the reference node in \autoref{fig:links} for demonstration purposes.
The respective model files that allow for results reproduction can be found in \cite{vlachas2021two}, as this example has been published as a benchmark multi-degree of freedom nonlinear response simulator.
Regarding material properties, the case study follows the template configuration \cite{vlachas_konstantinos_2021_4742248}.
Specifically, steel HEA cross-sections have been used for all beam elements, whereas the structure is assembled using two frames along axis $x$, each of $l = 7.5\:m$ length and one of $w = 6\:m$ along the width.
In addition, each story has a height of $h = 3.2\:m$.

A Bouc-Wen formulation has been utilised to model the behavior of the nonlinear joints: this reflects a smooth hysteretic model, often adopted for modeling material nonlinearity \cite{ikhouane2007systems}.
Therefore, based on the benchmark description in \cite{vlachas2021two} a Bouc-Wen model is introduced at every \ac{DOF} of every nodal coupling to simulate the total restoring force $\mathbf{R}$ of each joint.
An example illustration of the nonlinear mechanism in the longitudinal x-\ac{DOF} is provided in \autoref{fig:bw}.
The Bouc-Wen link models $\mathbf{R}$ as a superposition of a linear and a nonlinear term, represented by the two springs in \autoref{fig:bw}.
The linear and nonlinear terms, or springs, depend on the instantaneous nodal response $\delta u$ and the hysteretic, and thus history-dependent, component of the response $z$, respectively. 
In turn, the respective vectorized mathematical formulation for all \ac{DOF}s of the link reads:
\begin{align}
\mathbf{R} = \mathbf{R}_{linear} + \mathbf{R}_{hysteretic} =  \alpha*k*\mathbf{du} + (1-\alpha)*k*\mathbf{z} 
\label{eq:restforces}
\end{align}
where $d\mathbf{u}$ represents the nodal displacement, and $\alpha,k$ are traits characterizing the Bouc-Wen model on each link.
Regarding their physical interpretation, $\alpha$ represents the characteristic post-yield to elastic stiffness reaction for each link, whereas $k$ is the corresponding stiffness coefficient. 
The variable $\mathbf{z}$ stands for the hysteretic portion of the elongation, or displacement in general, and controls the hysteretic forcing. 
It obeys the following:
\begin{gather}
\dot{\mathbf{z}} = \frac{A \mathbf{d}\dot{\mathbf{u}} - \nu (t) (\beta |\mathbf{d}\dot{\mathbf{u}}|\mathbf{z}|\mathbf{z}|^{w-1}-\gamma \mathbf{d}\dot{\mathbf{u}}|\mathbf{z}|^{w})}{\eta (t)}  \label{eq:boucwendegrad}\\
\nu (t) = 1.0 + \delta_{\nu} \epsilon(t), \quad \eta (t) = 1.0 + \delta_{\eta} \epsilon(t), \quad \epsilon(t) = \int_0^{t} \mathbf{z}\mathbf{d}\dot{\mathbf{u}} dt 
\end{gather}
where the shape, smoothness, and overall amplitude of the hysteretic curve that characterises the dynamic behavior of each joint is determined by the Bouc-Wen parameters $\beta, \gamma, w$, and $A$ respectively.
The terms $\nu (t), \eta (t)$ are introduced to capture strength deterioration and stiffness degradation effects via the corresponding coefficients $\delta_{\nu}$ and $\delta_{eta}$.
In turn, their evolution in time depends on  the absorbed hysteretic energy, $\epsilon (t)$.
This representation allows for a structural dynamics simulator, which is parametrised with respect to system properties and traits of the joints' behavior.
For a more detailed elaboration on the physical connotations of the Bouc-Wen model parameters in terms of yielding, softening, and hysteretic behavior effects, the reader is referred to \cite{Chatzi2009, vlachas2021two}.

This parameterised shear frame simulator is selected due to its ability to model a variety of nonlinear dynamic effects that dominate the response and are dependent on the parametric configuration of the model.
In the presented case studies, the parameters defining the structure itself and those defining the acting loads can significantly affect the response.
The parameter set includes the forcing signal's temporal and spectral characteristics, the frame's material properties, and the traits that dictate the hysteretic effects on the joints.
The six parameters employed in this example are summarised in \autoref{tab:paramsBW}.

First, uncertainty is introduced in the material properties of the system by treating the Young modulus of elasticity $E$ as a parameter of the model.
Its range is summarized in \autoref{tab:paramsBW}.
In addition, three traits of the nonlinear joints of the shear frame are parametrically to model and simulate various qualities and shapes of the corresponding hysteretic behavioral curves.
Specifically, parameters $\alpha$ and $k$ in \autoref{eq:restforces} and parameter $\delta_n$ in \autoref{eq:boucwendegrad} are injected as dependencies in the derived \ac{ROM}.
The numerical range for each parameter is also provided in \autoref{tab:paramsBW}.

\begin{table}[!ht]
\caption{Two-story frame: Range of the parameter values of the \ac{ROM}. All parameters follow a uniform distribution. \label{tab:params}}
\label{tab:paramsBW}
\centering
\begin{tabular}{l c c c c c c c c }
\hline 
Parameter: & $\alpha$ & $k$ & $Amp$ & $f_{but}$ & $E$ & $\delta_{\eta}$ \\
Range: & [0.25,0.50] & [0.8,1.2]$\times 10^8$ & [1.5,3.0]$\times 10^6$ & [5,15] $Hz$ & [185,235] $GPa$ & [0.25, 0.75] \\
\hline
\end{tabular}
\end{table}

Forcing is applied to the frame system as a base excitation scenario representing an earthquake.
The force is applied at an angle of $\theta = \pi/4$ with respect to the x-axis, as depicted in \autoref{fig:links}.
To produce a parameterised version of the excitation, a white noise template accelerogram is used as a reference. 
This template signal is then passed through a second-order butterworth low-pass filter and multiplied with an amplitude factor to produce the actual accelerogram of the motion imposed on the system.
The amplitude coefficient $Amp$ and the frequency of the filter $f_{but}$ are treated as dependencies.
Thus, due to the dependencies of the model both in system parameters and excitation traits, the dynamic system under consideration exhibits substantially different behavior depending on the chosen parametric configuration.
The numerical study has been designed in such a way to validate the requirement to inject dependencies in the derived surrogate while making use of the ability of the simulator to output a variety of nonlinear behaviors, thus challenging the accuracy limit of the \ac{ROMs}. 


All \ac{ROM}s as referenced in \autoref{tab:pROMS} are implemented here, employing the same training scheme of fifty samples, drawn using \ac{LHS}.
The corresponding performance measures for each \ac{ROM} are evaluated on a validation set of five hundred samples, drawn using an \ac{LHS} with a different seed. 
Regarding the low-order dimension, $r=16$ modes are retained for each local basis $\bV$ and $\Tilde{r}=200$ modes for $\bV_{global}$ in \autoref{eq:coeffs}.

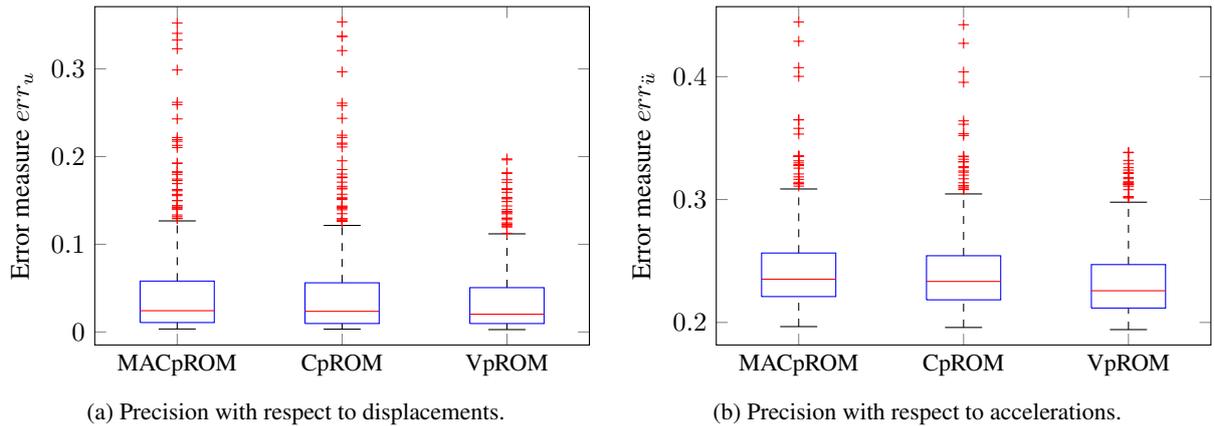
\begin{figure}[!ht]
    \centering
	\begin{subfigure}[c]{0.49\textwidth}
        \centering
		\begin{tikzpicture}

\begin{axis}[%
width=1.0\textwidth,
height=0.75\textwidth,
at = {(0,0)},
xmin=0.5,
xmax=3.5,
xtick={1,2,3},
xticklabels={{MACpROM},{CpROM},{VpROM}},
xticklabel style = {font =\fontsize{\figureFontSize pt}{9pt}\selectfont},
ymin=-0.0148545021139751,
ymax=0.371030659909914,
ylabel={Error measure $err_u$},
scaled ticks=false,
axis background/.style={fill=white}
]
\addplot [color=black, dashed, forget plot]
  table[row sep=crcr]{%
1	0.0580338271080829\\
1	0.126641900852208\\
};
\addplot [color=black, dashed, forget plot]
  table[row sep=crcr]{%
2	0.0560478125764086\\
2	0.121497005415242\\
};
\addplot [color=black, dashed, forget plot]
  table[row sep=crcr]{%
3	0.0505230486072289\\
3	0.111889493750408\\
};
\addplot [color=black, dashed, forget plot]
  table[row sep=crcr]{%
1	0.00339560207489808\\
1	0.0107355935014033\\
};
\addplot [color=black, dashed, forget plot]
  table[row sep=crcr]{%
2	0.00323168148484029\\
2	0.00963088213209778\\
};
\addplot [color=black, dashed, forget plot]
  table[row sep=crcr]{%
3	0.00268573252347439\\
3	0.00958181843475263\\
};
\addplot [color=black, forget plot]
  table[row sep=crcr]{%
0.8875	0.126641900852208\\
1.1125	0.126641900852208\\
};
\addplot [color=black, forget plot]
  table[row sep=crcr]{%
1.8875	0.121497005415242\\
2.1125	0.121497005415242\\
};
\addplot [color=black, forget plot]
  table[row sep=crcr]{%
2.8875	0.111889493750408\\
3.1125	0.111889493750408\\
};
\addplot [color=black, forget plot]
  table[row sep=crcr]{%
0.8875	0.00339560207489808\\
1.1125	0.00339560207489808\\
};
\addplot [color=black, forget plot]
  table[row sep=crcr]{%
1.8875	0.00323168148484029\\
2.1125	0.00323168148484029\\
};
\addplot [color=black, forget plot]
  table[row sep=crcr]{%
2.8875	0.00268573252347439\\
3.1125	0.00268573252347439\\
};
\addplot [color=blue, forget plot]
  table[row sep=crcr]{%
0.775	0.0107355935014033\\
0.775	0.0580338271080829\\
1.225	0.0580338271080829\\
1.225	0.0107355935014033\\
0.775	0.0107355935014033\\
};
\addplot [color=blue, forget plot]
  table[row sep=crcr]{%
1.775	0.00963088213209778\\
1.775	0.0560478125764086\\
2.225	0.0560478125764086\\
2.225	0.00963088213209778\\
1.775	0.00963088213209778\\
};
\addplot [color=blue, forget plot]
  table[row sep=crcr]{%
2.775	0.00958181843475263\\
2.775	0.0505230486072289\\
3.225	0.0505230486072289\\
3.225	0.00958181843475263\\
2.775	0.00958181843475263\\
};
\addplot [color=red, forget plot]
  table[row sep=crcr]{%
0.775	0.024223035195391\\
1.225	0.024223035195391\\
};
\addplot [color=red, forget plot]
  table[row sep=crcr]{%
1.775	0.0235122317153441\\
2.225	0.0235122317153441\\
};
\addplot [color=red, forget plot]
  table[row sep=crcr]{%
2.775	0.0202635508824856\\
3.225	0.0202635508824856\\
};
\addplot [color=black, only marks, mark=+, mark options={solid, draw=red}, forget plot]
  table[row sep=crcr]{%
1	0.12928304814547\\
1	0.129595644195053\\
1	0.131461654533375\\
1	0.133418314871538\\
1	0.139880366620071\\
1	0.141163462858978\\
1	0.142109062713807\\
1	0.14390708758054\\
1	0.143941114370423\\
1	0.149767771709564\\
1	0.149856465376588\\
1	0.155461514204814\\
1	0.156900143927019\\
1	0.161065505992843\\
1	0.162461385181278\\
1	0.169201272429875\\
1	0.172731818161221\\
1	0.174372380503657\\
1	0.179750682820698\\
1	0.182075307538043\\
1	0.182687449567239\\
1	0.192235476604795\\
1	0.192637254215588\\
1	0.21033671594931\\
1	0.212748676767009\\
1	0.217425383372236\\
1	0.219288387288175\\
1	0.221725615360165\\
1	0.242981297425415\\
1	0.259251346581173\\
1	0.261895471768105\\
1	0.298806978468088\\
1	0.322961713717751\\
1	0.332918368619655\\
1	0.340548168923564\\
1	0.352295800849501\\
};
\addplot [color=black, only marks, mark=+, mark options={solid, draw=red}, forget plot]
  table[row sep=crcr]{%
2	0.125975682611609\\
2	0.126197584267377\\
2	0.126590756344787\\
2	0.127483444886914\\
2	0.128937979154929\\
2	0.129195747545042\\
2	0.134857411616545\\
2	0.139328887835\\
2	0.141039026508172\\
2	0.142163068545086\\
2	0.142715942473421\\
2	0.1436364948824\\
2	0.15101908493386\\
2	0.152438349211589\\
2	0.152994359032874\\
2	0.157056094240635\\
2	0.163103428884516\\
2	0.164008193107256\\
2	0.166100541986473\\
2	0.170859925081604\\
2	0.175943382326531\\
2	0.176191645861007\\
2	0.176613766554718\\
2	0.180040107647042\\
2	0.185229717572643\\
2	0.195186555284902\\
2	0.210943665473701\\
2	0.214176794095368\\
2	0.21534622220615\\
2	0.221758479387685\\
2	0.224379270999533\\
2	0.243649842256233\\
2	0.25818267106275\\
2	0.260991928849603\\
2	0.296610213307725\\
2	0.320873902152221\\
2	0.337008036259374\\
2	0.337250118449526\\
2	0.353490425272464\\
};
\addplot [color=black, only marks, mark=+, mark options={solid, draw=red}, forget plot]
  table[row sep=crcr]{%
3	0.11255078738183\\
3	0.119784432152162\\
3	0.12009007417223\\
3	0.121848371283339\\
3	0.12260096420397\\
3	0.123113517214435\\
3	0.124408651816542\\
3	0.128043508072842\\
3	0.129531857581281\\
3	0.12990369162012\\
3	0.13476216702362\\
3	0.135962883448072\\
3	0.137088568147071\\
3	0.137384116732235\\
3	0.139584838915175\\
3	0.143693481303689\\
3	0.148706577787396\\
3	0.151845793464614\\
3	0.152885979860907\\
3	0.153365224013001\\
3	0.159249156333021\\
3	0.162385269904048\\
3	0.163075788234279\\
3	0.163295561378207\\
3	0.170617548343885\\
3	0.173616202964191\\
3	0.180791242469538\\
3	0.181610067338936\\
3	0.196287188255006\\
3	0.197422004604454\\
};
\end{axis}
\end{tikzpicture}%

		\caption{Precision with respect to displacements.}
		\label{fig:AccDisps}
	\end{subfigure}
	\begin{subfigure}[c]{0.49\textwidth}
        \centering
		\begin{tikzpicture}

\begin{axis}[%
width=1.0\textwidth,
height=0.75\textwidth,
at = {(0,0)},
xmin=0.5,
xmax=3.5,
xtick={1,2,3},
xticklabels={{MACpROM},{CpROM},{VpROM}},
ymin=0.181675964595168,
ymax=0.456941380967836,
xticklabel style = {font =\fontsize{\figureFontSize pt}{9pt}\selectfont},
ylabel={Error measure $err_{\ddot{u}}$},
axis background/.style={fill=white}
]
\addplot [color=black, dashed, forget plot]
  table[row sep=crcr]{%
1	0.256481887411275\\
1	0.30859729834966\\
};
\addplot [color=black, dashed, forget plot]
  table[row sep=crcr]{%
2	0.254285720026822\\
2	0.304563998830478\\
};
\addplot [color=black, dashed, forget plot]
  table[row sep=crcr]{%
3	0.247084590924984\\
3	0.297843751136304\\
};
\addplot [color=black, dashed, forget plot]
  table[row sep=crcr]{%
1	0.196628247145831\\
1	0.221053471226956\\
};
\addplot [color=black, dashed, forget plot]
  table[row sep=crcr]{%
2	0.195906691193779\\
2	0.218400648191911\\
};
\addplot [color=black, dashed, forget plot]
  table[row sep=crcr]{%
3	0.194188028975744\\
3	0.211729811545903\\
};
\addplot [color=black, forget plot]
  table[row sep=crcr]{%
0.8875	0.30859729834966\\
1.1125	0.30859729834966\\
};
\addplot [color=black, forget plot]
  table[row sep=crcr]{%
1.8875	0.304563998830478\\
2.1125	0.304563998830478\\
};
\addplot [color=black, forget plot]
  table[row sep=crcr]{%
2.8875	0.297843751136304\\
3.1125	0.297843751136304\\
};
\addplot [color=black, forget plot]
  table[row sep=crcr]{%
0.8875	0.196628247145831\\
1.1125	0.196628247145831\\
};
\addplot [color=black, forget plot]
  table[row sep=crcr]{%
1.8875	0.195906691193779\\
2.1125	0.195906691193779\\
};
\addplot [color=black, forget plot]
  table[row sep=crcr]{%
2.8875	0.194188028975744\\
3.1125	0.194188028975744\\
};
\addplot [color=blue, forget plot]
  table[row sep=crcr]{%
0.775	0.221053471226956\\
0.775	0.256481887411275\\
1.225	0.256481887411275\\
1.225	0.221053471226956\\
0.775	0.221053471226956\\
};
\addplot [color=blue, forget plot]
  table[row sep=crcr]{%
1.775	0.218400648191911\\
1.775	0.254285720026822\\
2.225	0.254285720026822\\
2.225	0.218400648191911\\
1.775	0.218400648191911\\
};
\addplot [color=blue, forget plot]
  table[row sep=crcr]{%
2.775	0.211729811545903\\
2.775	0.247084590924984\\
3.225	0.247084590924984\\
3.225	0.211729811545903\\
2.775	0.211729811545903\\
};
\addplot [color=red, forget plot]
  table[row sep=crcr]{%
0.775	0.235205094220746\\
1.225	0.235205094220746\\
};
\addplot [color=red, forget plot]
  table[row sep=crcr]{%
1.775	0.233355281402352\\
2.225	0.233355281402352\\
};
\addplot [color=red, forget plot]
  table[row sep=crcr]{%
2.775	0.225777934210061\\
3.225	0.225777934210061\\
};
\addplot [color=black, only marks, mark=+, mark options={solid, draw=red}, forget plot]
  table[row sep=crcr]{%
1	0.310909696845153\\
1	0.312812239690815\\
1	0.312952150888619\\
1	0.314016665177168\\
1	0.316386924508178\\
1	0.318448952721896\\
1	0.320823387305765\\
1	0.325533178844879\\
1	0.327668670733812\\
1	0.328420850234845\\
1	0.329817922686961\\
1	0.33151510914371\\
1	0.334966796773361\\
1	0.334968994583401\\
1	0.335282161072934\\
1	0.335600574494445\\
1	0.353274015428817\\
1	0.357913885812751\\
1	0.364789086767965\\
1	0.36501230291796\\
1	0.400204394684659\\
1	0.407280031205197\\
1	0.428775160858694\\
1	0.444429316587261\\
};
\addplot [color=black, only marks, mark=+, mark options={solid, draw=red}, forget plot]
  table[row sep=crcr]{%
2	0.308119407209033\\
2	0.308411434102113\\
2	0.309567877216882\\
2	0.311059863827504\\
2	0.311275065693454\\
2	0.313653196474418\\
2	0.31687524477167\\
2	0.319892229574432\\
2	0.322207370552472\\
2	0.323273801500815\\
2	0.325916753673463\\
2	0.326873624288522\\
2	0.330620180498039\\
2	0.330893531514777\\
2	0.332775716189505\\
2	0.334860734993762\\
2	0.335288991262681\\
2	0.352084539616993\\
2	0.353520216749418\\
2	0.361255379926535\\
2	0.364015371029168\\
2	0.39538588033965\\
2	0.403803584212582\\
2	0.427102062587802\\
2	0.442213645111549\\
};
\addplot [color=black, only marks, mark=+, mark options={solid, draw=red}, forget plot]
  table[row sep=crcr]{%
3	0.301249603218319\\
3	0.301868788876303\\
3	0.302480698797575\\
3	0.308027896955497\\
3	0.310840066182755\\
3	0.312025112437845\\
3	0.312407479667713\\
3	0.312440155028996\\
3	0.313198960769161\\
3	0.313691722911037\\
3	0.314502865018969\\
3	0.316952397492609\\
3	0.317947349316427\\
3	0.321156705101873\\
3	0.321988618103567\\
3	0.32396652478874\\
3	0.32593386921976\\
3	0.329036823826835\\
3	0.331844047008191\\
3	0.338144309086923\\
3	0.338459449111884\\
};
\end{axis}
\end{tikzpicture}%
		\caption{Precision with respect to accelerations.}
		\label{fig:AccAs}
	\end{subfigure}
 \caption{Box plots reporting the accuracy of capturing the displacement and acceleration time histories. The distributions of the respective error measures $err_u,err_{\ddot{u}}$ from \autoref{eq:errors} are visualised along with the respective median (red line) and outliers (red crosses). }
	\label{fig:BoxPlots}
\end{figure}
\noindent
A detailed evaluation of the accuracy for the implemented \ac{ROM}s of \autoref{tab:pROMS} is presented in \autoref{fig:BoxPlots}.
The precision of the respective surrogates is evaluated with respect to two measures, namely $err_u$ and $err_{\ddot{u}}$ of \autoref{eq:errors} that correspond to the error on capturing the displacement and acceleration time histories respectively. 
The boxplots provide a visualisation of the ability of each \ac{ROM} to capture the \ac{FOM} response both in terms of displacements and accelerations, whereas the respective values are also reported in \autoref{tab:ErrorsBW}.
Although the overall precision is relatively low,
this example merely serves to offer a comparison that deliberately employs a relatively wide domain of parameters, in order to excite substantially different dynamic behavior.

\begin{table}[!htb]
    \centering
	\caption{Performance for the \ac{ROM}s from \autoref{tab:pROMS}. The median and maximum error is reported with respect to displacements and accelerations. Efficiency is also reported, although hyper-reduction has not been implemented.}
	\label{tab:ErrorsBW}
		\begin{tabular} { c|c c c c|c c| }
        \hline
        & \multicolumn{2}{c}{Median}
        & \multicolumn{2}{c}{Maximum}
        & \multicolumn{2}{c}{Efficiency}
        \\
        & \multicolumn{2}{c}{Error}
        & \multicolumn{2}{c}{Error}
        & Speed-Up
        & CPU
        \\
        & $err_u$ & $err_{\ddot{u}}$ & $err_u$ & $err_{\ddot{u}}$ & Factor & timing
        \\
		\hline
        \ac{FOM} & - & - & - & - & 1.00 & 53 (s) \\
		\ac{MACpROM} & 2.42$\%$ & 35.23$\%$ & 23.52$\%$ & 44.44$\%$ & 2.22 & 24 (s) \\
		CpROM & 2.35$\%$ & 35.34$\%$ & 23.33$\%$ & 44.22$\%$ & 2.10 & 25 (s) \\
		VpROM & 2.20$\%$ & 19.62$\%$ & 22.57$\%$ & 33.81$\%$ & 2.41 & 22 (s) \\
	    \hline
    \end{tabular}
\end{table}

Nevertheless, in \autoref{fig:BoxPlots} and \autoref{tab:ErrorsBW} the \ac{MACpROM} implemented exhibits similar accuracy with the reference CpROM, in terms of approximating both the displacement and the acceleration time histories. 
The respective median error and the boxplot quartiles almost coincide, whilst both approaches seem to deliver a similar distribution of outliers. 
The proposed \emph{VpROM} on the other hand achieves a substantially improved performance.
The outliers are fewer, the respective discrepancy for the outlier samples is substantially lower, and the visualised distribution has a lower median and maximum error. 

\begin{figure}[!hb]
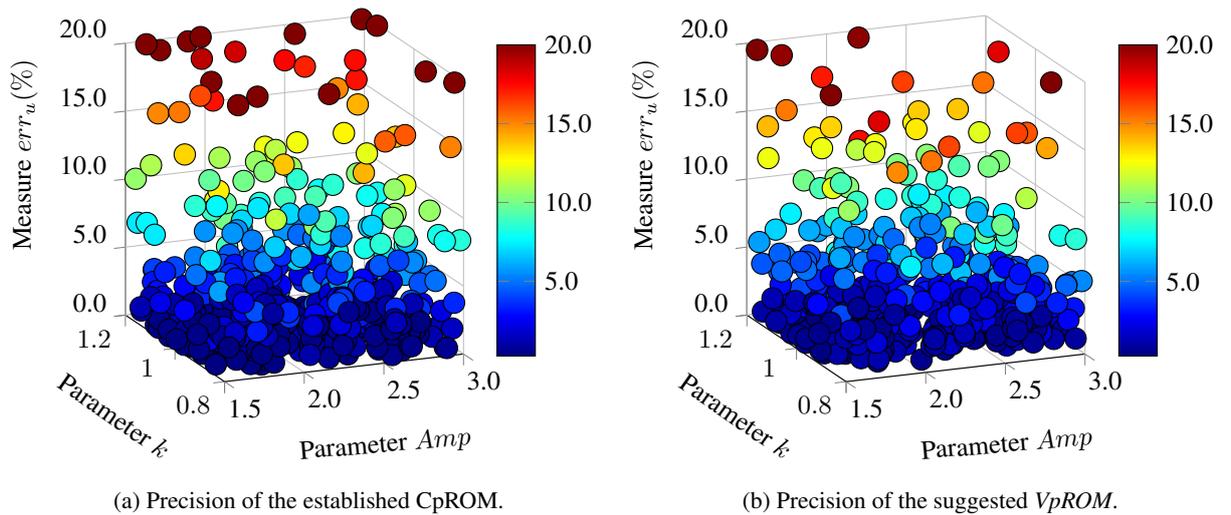

    \centering
	\begin{subfigure}[b]{0.49\textwidth}
        \input{figures_tikz/Figure7a_ErrorsCpROM.tikz}
		\caption{Precision of the established CpROM.}
		\label{fig:ErrorsMAC}
	\end{subfigure}
	\begin{subfigure}[b]{0.49\textwidth}
		\input{figures_tikz/Figure7b_ErrorsVpROM.tikz}
        \caption{Precision of the suggested \emph{VpROM}.}
		\label{fig:ErrorsVAE}
	\end{subfigure}
	\caption{Parameter Visualisation of the error distribution of \autoref{eq:errors} for displacement time histories ($err_u$). The proposed \emph{VpROM} is compared with CpROM, as the MACpROM delivers a slightly worse performance as indicated in \autoref{fig:BoxPlots}. All samples with errors greater than $20\%$ are depicted in the $20\%$ color level for a clearer comparison.}	\label{fig:Errors}
\end{figure}

A further comparative visualisation of the accuracy for the implemented \ac{ROM}s is provided in \autoref{fig:Errors}.
An example projection plane has been chosen for demonstration purposes and the validation measure is depicted on the vertical axis and via the color scale.
Similar to \autoref{fig:AccDisps}, the $err_u$ is visualised in \autoref{fig:Errors} as a representative measure of the \ac{ROM}s ability to reproduce the displacement time histories of the \ac{FOM}.
In addition, all samples with errors greater than $20\%$ are depicted in the $20\%$ color level for better scaling and a clearer comparison. 
Since the established CpROM and the \ac{MACpROM} deliver similar precision with respect to displacements in \autoref{fig:AccDisps}, the \emph{VpROM} suggested in this study is compared only with CpROM in \autoref{fig:Errors} for the sake of a clearer demonstration.

As already highlighted, the \emph{VpROM} captures the dynamics across the domain of parametric inputs with an overall superior precision and fewer accuracy outliers than CpROM or \ac{MACpROM}.
This is visualised in \autoref{fig:Errors} through the fewer circles located in the dark red region for the \emph{VpROM} and the substantially fewer evaluations colored outside the blue-to-green spectrum.
Thus, despite a few outliers, the overall accuracy of the framework remains superior to the compared established alternatives for physics-based \ac{MOR}. 

In \autoref{fig:RespOverview}, a more detailed evaluation of the approximation quality achieved by the \emph{VpROM} is illustrated.
Specifically, the time history estimation is depicted for different levels of precision to visualise and validate the overall performance of the \emph{VpROM}.
One sample from each family of sample points as captured by the color scale in \autoref{fig:Errors} is presented. 
The \emph{VpROM} approximation is visualised for various response patterns to highlight the ability of the surrogate model to infer dynamic behaviors dominated by different effects as modeled via the variety of shapes of the hysteretic curves on the links.
The \emph{VpROM} is shown to deliver a robust approximation in a complex example with a rich dynamic behavior represented by several different shapes and amplitudes of the hysteretic curves characterizing the behavior of the nonlinear joints.

\begin{figure}[!htb]
\centering
    \begin{subfigure}[c]{0.49\textwidth}
        \centering
        \begin{tikzpicture}
	\begin{axis}[
		name = response,
		xmin = 0,
		xmax = 10.0,
		ymin = -0.4,
		ymax = 0.4,
		xtick = {0, 2, 4, 6, 8,10},
		ytick = {-0.4,-0.2, 0, 0.2, 0.4},
		xlabel = {Time (s)},
		ylabel = {Displacement ($m$)},
		grid = both,
		width=\textwidth,
		height=0.55\textwidth,
        legend style={legend columns=-1, legend pos=south east,nodes={font=\fontsize{\figureFontSize pt}{\figureFontSize pt}\selectfont}}
		]
			
		\addplot[color = blue , line width=0.75pt] table[x=time, y=one] {figures_tikz/THdata.dat};

        \addplot[color = red, dashed, line width=0.75pt] table[x=time, y=seven] {figures_tikz/THdata.dat};
  		
		\addlegendentry{FOM}
		\addlegendentry{VpROM}
	\end{axis}
	
\end{tikzpicture}
        \captionsetup{skip=-8pt}
        \caption{Approximation for a sample in the dark red region in Fig. \ref{fig:Errors}.}  \label{fig:sampleA}
    \end{subfigure}     
    \begin{subfigure}[c]{0.49\textwidth}
        \centering
        \begin{tikzpicture}
	\begin{axis}[
		name = response,
		xmin = 0,
		xmax = 10.0,
		ymin = -0.5,
		ymax = 0.5,
		xtick = {0, 2, 4, 6, 8,10},
		ytick = {-0.5,-0.25, 0, 0.25, 0.5},
		xlabel = {Time (s)},
		ylabel = {Displacement ($m$)},
		grid = both,
		width=\textwidth,
		height=0.55\textwidth,
        legend style={legend columns=-1, legend pos=south east,nodes={font=\fontsize{\figureFontSize pt}{\figureFontSize pt}\selectfont}}
		]
			
		\addplot[color = blue , line width=0.75pt] table[x=time, y=four] {figures_tikz/THdata.dat};

        \addplot[color = red, dashed, line width=0.75pt] table[x=time, y=ten] {figures_tikz/THdata.dat};
  		
		\addlegendentry{FOM}
		\addlegendentry{VpROM}

	\end{axis}
	
\end{tikzpicture}
        \captionsetup{skip=-8pt}
        \caption{Approximation for a sample in the red region in \autoref{fig:Errors}.}\label{fig:sampleB}
    \end{subfigure}
     \begin{subfigure}[c]{0.49\textwidth}
        \centering
        \begin{tikzpicture}
	\begin{axis}[
		name = response,
		xmin = 0,
		xmax = 10.0,
		ymin = -0.4,
		ymax = 0.4,
		xtick = {0, 2, 4, 6, 8,10},
		ytick = {-0.4,-0.2, 0, 0.2, 0.4},
		xlabel = {Time (s)},
		ylabel = {Displacement ($m$)},
		grid = both,
		width=\textwidth,
		height=0.55\textwidth,
		legend style={legend columns=-1, legend pos=south west,nodes={font=\fontsize{\figureFontSize pt}{\figureFontSize pt}\selectfont}}
		]
			
		\addplot[color = blue , line width=0.75pt] table[x=time, y=three] {figures_tikz/THdata.dat};

        \addplot[color = red, dashed, line width=0.75pt] table[x=time, y=nine] {figures_tikz/THdata.dat};
  		
		\addlegendentry{FOM}
		\addlegendentry{VpROM}

	\end{axis}
	
\end{tikzpicture}
        \captionsetup{skip=-8pt}
        \caption{Approximation for a sample in the yellow region in \autoref{fig:Errors}.} \label{fig:sampleC}
    \end{subfigure}     
    \begin{subfigure}[c]{0.49\textwidth}
        \centering
        \begin{tikzpicture}
	\begin{axis}[
		name = response,
		xmin = 0,
		xmax = 10.0,
		ymin = -0.6,
		ymax = 0.6,
		xtick = {0, 2, 4, 6, 8,10},
		ytick = {-0.6,-0.3, 0, 0.3, 0.6},
		xlabel = {Time (s)},
		ylabel = {Displacement ($m$)},
		grid = both,
		width=\textwidth,
		height=0.55\textwidth,
        legend style={legend columns=-1, legend pos=south west,nodes={font=\fontsize{\figureFontSize pt}{\figureFontSize pt}\selectfont}}
		]
			
		\addplot[color = blue , line width=0.75pt] table[x=time, y=six] {figures_tikz/THdata.dat};

        \addplot[color = red, dashed, line width=0.75pt] table[x=time, y=twelve] {figures_tikz/THdata.dat};
  		
		\addlegendentry{FOM}
		\addlegendentry{VpROM}

	\end{axis}
	
\end{tikzpicture}
        \captionsetup{skip=-8pt}
        \caption{Approximation for a sample in the cyan region in \autoref{fig:Errors}.} \label{fig:sampleD}
    \end{subfigure}
	\caption{Visualisation of the different levels of the approximation quality achieved using the \emph{VpROM}. The \emph{VpROM} estimation is reported for various response patterns the system exhibits depending on its parametric features.}
	\label{fig:RespOverview}
\end{figure}
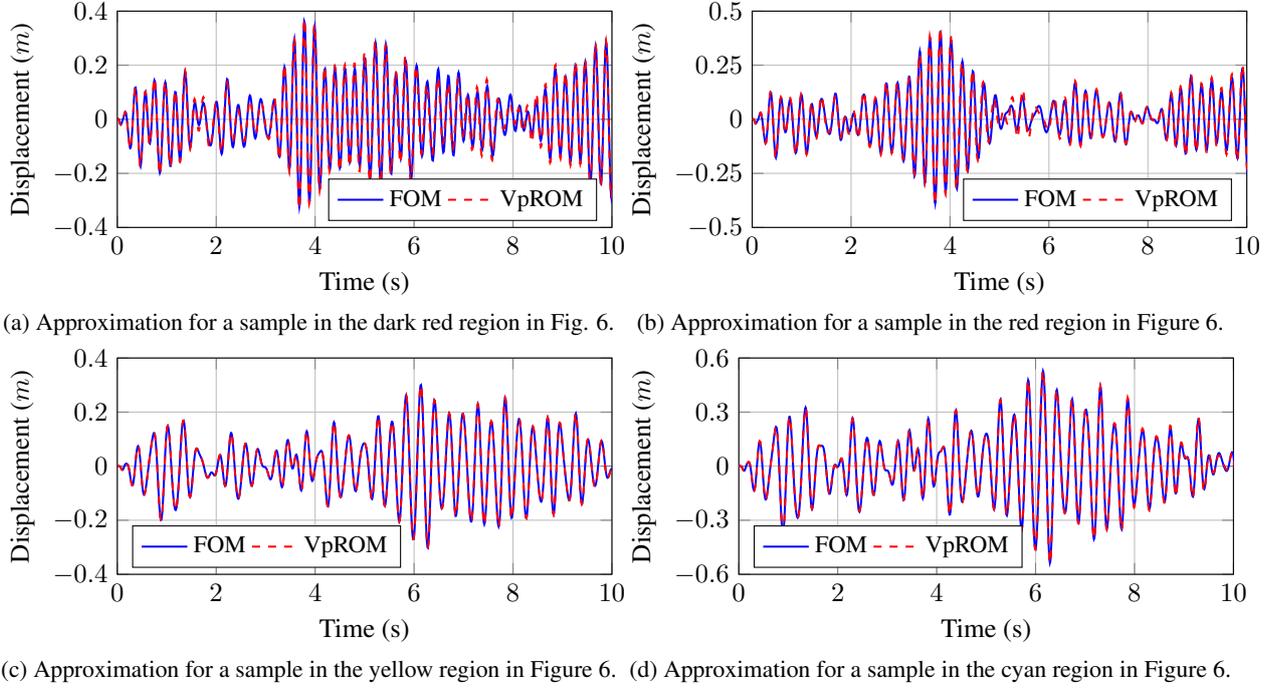


\begin{figure}[!b]
	\centering
 \begin{subfigure}[b]{0.49\textwidth}
 \centering
 \begin{tikzpicture}[spy using outlines=
	{rectangle, magnification=2, anchor = center, connect spies}]
	\begin{axis}[
		name = response,
		xmin = 3,
		xmax = 9,
		ymin = -0.25,
		ymax = 0.25,
		xtick = {0,3,5,7,9,10},
		ytick = {-0.2,0.0,0.2},
		xlabel = {Time (s)},
		ylabel = {Displacement ($m$)},
		grid = both,
		width=1.00\textwidth,
		height=0.60\textwidth,
		legend style={legend columns=-1, legend pos=south east,nodes={font=\fontsize{\figureFontSize pt}{\figureFontSize pt}\selectfont}}
		]
		
		\addplot[name path=F, color = red, line width=0.55pt] table[x=Time, y=FOM] {figures_tikz/UQplot.dat};
		\addplot[color = blue, dashed, line width=0.55pt] table[x=Time, y=ROM] {figures_tikz/UQplot.dat};
		\addplot[name path=G, color = blue,opacity=0.2] table[x=Time, y=Maxes] {figures_tikz/UQplot.dat};			
		\tikzfillbetween[of=F and G]{blue, opacity=0.2};
		\addplot[name path=G, color = blue,opacity=0.2] table[x=Time, y=Infs] {figures_tikz/UQplot.dat};			
		\tikzfillbetween[of=F and G]{blue, opacity=0.2};
		\addlegendentry{FOM}
		\addlegendentry{VpROM}

		
	\end{axis}
	
	
\end{tikzpicture}
\captionsetup{skip=4pt}
\caption{Approximation on displacement response.}
 \label{fig:THDisps}
 \end{subfigure}
 \begin{subfigure}[b]{0.49\textwidth}
 \centering
\begin{tikzpicture}[spy using outlines=
	{rectangle, magnification=2, anchor = center, connect spies}]
	\begin{axis}[
		name = response,
		xmin = 3,
		xmax = 9,
		ymin = -6,
		ymax = 6,
		xtick = {0,3,5,7,9},
		ytick = {-6.0,0.0,6.0},
		xlabel = {Time (s)},
		ylabel = {Acceleration ($m/sec^2$)},
		grid = both,
		width=1.00\textwidth,
		height=0.60\textwidth,
		legend style={legend columns=-1, legend pos=south east,nodes={font=\fontsize{\figureFontSize pt}{\figureFontSize pt}\selectfont}}
		]
		
		\addplot[name path=F, color = black, line width=0.55pt] table[x=Time, y=FOM] {figures_tikz/UQplotA.dat};
		\addplot[color = orange, dashed, line width=0.55pt] table[x=Time, y=ROM] {figures_tikz/UQplotA.dat};
		\addplot[name path=G, color = green,opacity=0.2] table[x=Time, y=Maxes] {figures_tikz/UQplotA.dat};			
		\tikzfillbetween[of=F and G]{green, opacity=0.2};
		\addplot[name path=G, color = green,opacity=0.2] table[x=Time, y=Infs] {figures_tikz/UQplotA.dat};			
		\tikzfillbetween[of=F and G]{green, opacity=0.2};
		\addlegendentry{FOM}
		\addlegendentry{VpROM}

		
	\end{axis}
	
	
\end{tikzpicture}
 \caption{Approximation on acceleration response.}
 \label{fig:THacceleration}
 \end{subfigure}
  \caption{Average quality of the \emph{VpROM} approximation. Evaluation is performed on the degree of freedom with the maximum absolute response. The shaded area quantifies the uncertainty of the response inference.}
 \label{fig:THs}
\end{figure}
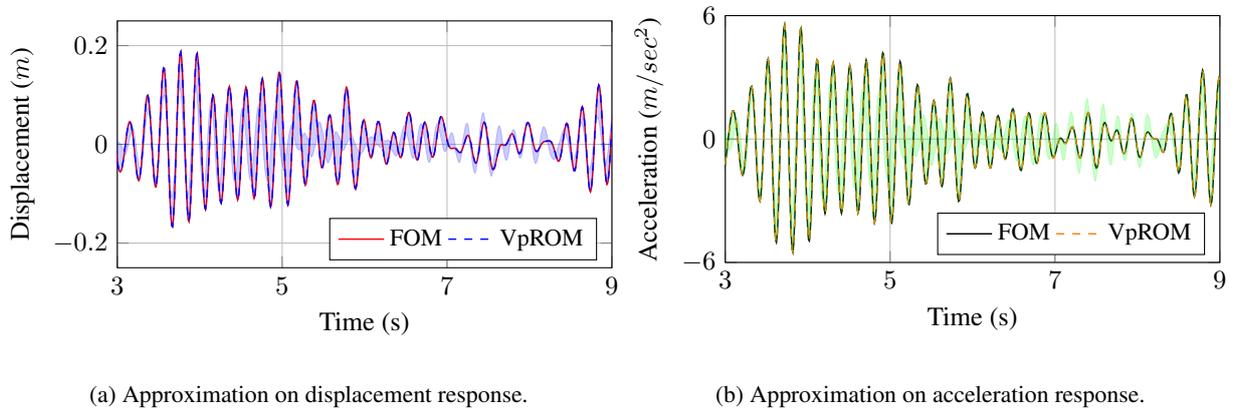

Beyond offering a robust reduction framework, that is shown to generalise across parametric configurations, the derived \emph{VpROM} offers the potential of quantifying the uncertainty in the respective estimations.
This is due to the latent space of the \ac{cVAE} component being trained to approximate a given variational distribution; thus, by sampling this variational distribution, the uncertainty on the predicted local basis coefficients may be captured in addition to simply mean estimates. 
Naturally, by using only the mean predicted values of the coefficients, significant information is lost regarding the uncertainty of these values.
Hence, it is worth considering utilising this technique for uncertainty estimation on the response of the \ac{ROM} for each sample.
In order to propagate this uncertainty, multiple parallel evaluations of the \emph{VpROM} are performed employing different coefficient values as generated from different samples from the latent space of the \ac{VAE}.
In turn, the distribution of the responses is inspected to evaluate and quantify the uncertainty in the respective inference.
The resulting approximation is equipped with confidence bounds that can provide increased utility for many problems in structural dynamics.

A visualisation of the respective output is provided in \autoref{fig:THs}, where the average performance of the \emph{VpROM} for both displacements and accelerations is depicted.
The respective shaded area represents the confidence bounds of the inference scheme, evaluated by sampling the predicted distributions of the latent space 40 times and propagating the response using the respective local bases assembled by the decoder for each of these 40 sampled vectors. The shaded area encompasses the maximum and minimum values at each time point for the 40 simulations carried out.
The respective average quality of the \emph{VpROM} approximation in \autoref{fig:THs} indicates a high-precision physics-based surrogate, with the inherent ability to provide a quantification on the uncertainty of the respective estimations.
To further evaluate the suitability of the proposed framework and demonstrate its utility in reducing computational toll, a large-scale example is discussed next.

\subsection{Wind Turbine Tower with plasticity}

\begin{figure}[!hb]
	\begin{subfigure}[b]{0.5\textwidth}
		\centering
		\includegraphics[scale=0.25]{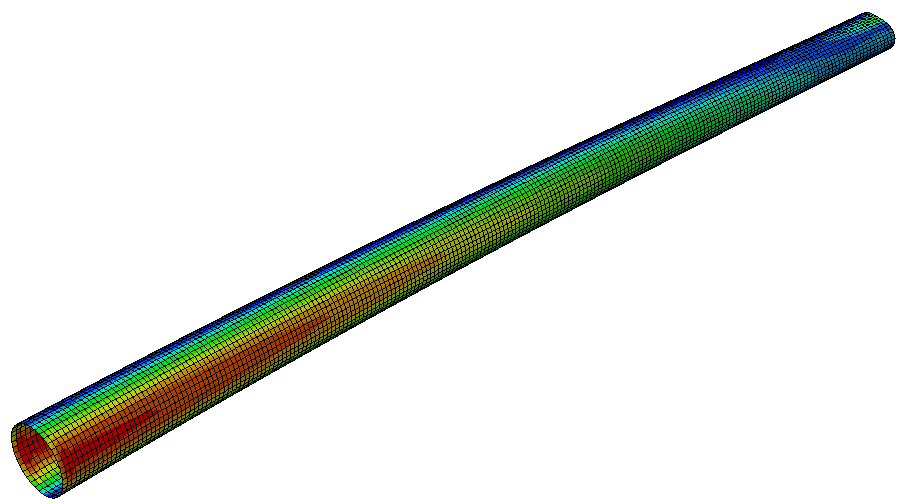}
		\caption{\ac{FE} mesh in an example deformed state.}
		\label{fig:geometryF}
	\end{subfigure}
    \begin{subfigure}[b]{0.5\textwidth}
		\centering
		\includegraphics[scale=0.25]{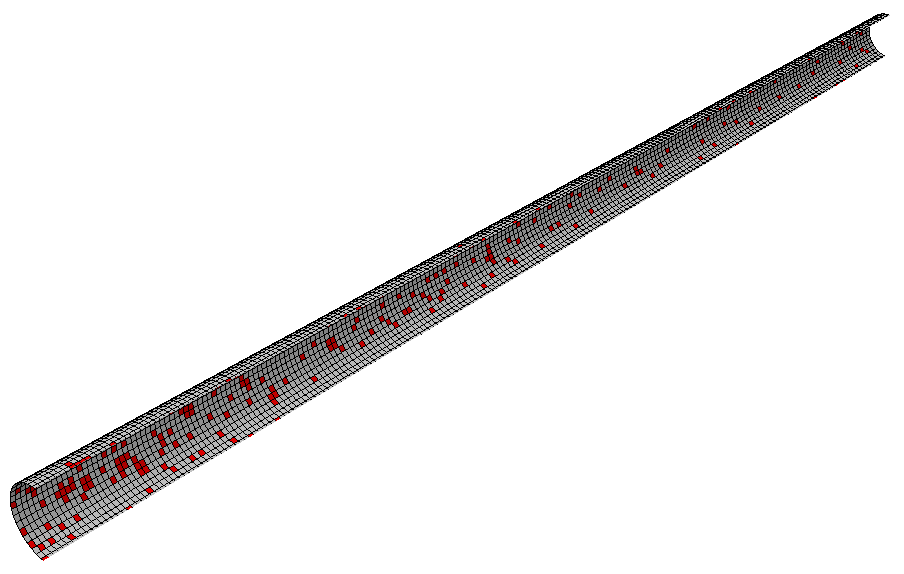}
		\caption{ECSW elements highlighted in red for one of the clusters.}
		\label{fig:hyperF}
	\end{subfigure}
	\caption{Wind turbine tower: \ac{FE} model and example of ECSW mesh. For the ECSW elements a horizontal cut is depicted for visualization purposes.}   
    \label{fig:geometry}
\end{figure}

This section evaluates the performance of the suggested \emph{VpROM} on a large-scale example based on the simulated dynamic response of the NREL 5-MW reference wind turbine tower \cite{Jonkman2009}.
Regarding the configuration of this case study, the interested reader is referred to \cite{Vlachas2021}.
In brief, the three-dimensional \ac{FE} model of the monopile is visualized in \autoref{fig:geometryF} and features a circular cross-section, which is linearly tapered from the base to the top. 
The respective diameter and wall thickness are equal to $6m$ and $0.027m$ on the base and $3.87m$ and $0.019m$ on the top of the monopile.
However, for simplification purposes, a constant thickness assumption is made throughout the tower, and $8170$ shell elements are used.
The wind turbine is assembled at the top of the monopile assuming a lumped mass scheme and regular beam elements, assembled through multi-point constraints to the tower. 
Regarding the material properties, steel is assumed, with $E_{steel}=210 GPa$ and a density of $\rho=7850 kg/m^3$ and a nonlinear constitutive law, which is characterised by isotropic von Mises plasticity. 

Although the employed stress-strain relation might seem relatively simple, this problem features an extensive yielding domain ($\approx 30\%$ of the height), additional model uncertainties, and a stochastic excitation that increase the complexity and pose certain requirements when deriving a high precision \ac{ROM}.
In addition, this large-scale case study is used to demonstrate the efficiency of hyper-reduced \ac{ROM}s, which allow for a substantial reduction in the overall computational toll.
For this reason, the hyper-reduced variant of all \ac{ROM}s, termed as \emph{HR-VpROM} for example in \autoref{tab:pROMS}, is validated herein.

\begin{table}[!ht]
\caption{Range of the parametric values of the implemented \ac{ROM}s. The $E_{ref}$ refers to the typical modulus assigned to each material of the model.}
\label{tab:paramsMonopile}
\centering
\begin{tabular}{l c c c }
\hline 
Parameter: & Excitation amplitude $A$ & Yield stress $\sigma_{VM}$ & Young modulus $E$ \\
Range: & [2.50, 3.75] & [375, 450] ($MPa$) & [0.80, 1.20]$\times E_{ref}$ ($GPa$) \\
\hline
\end{tabular}
\end{table}

Regarding the parametric dependencies, the Kobe earthquake accelerogram is utilised as a ground motion scenario, parameterised with respect to its amplitude $A$.
The yield stress $\sigma_{VM}$ and the Young modulus of elasticity $E$ are also varied.
The range of these parametric dependencies is summarised in \autoref{tab:paramsMonopile}.
The training and validation domain is designed using \ac{LHS} sampling, similar to \autoref{example:BWlinks}.
In this case study, a low-order dimension of $r=4$ is chosen, while $\Tilde{r}=32$ global modes are retained for $\bV_{global}$
in \autoref{eq:coeffs}, whereas the $\tau$ parameter for the ECSW hyper-reduction technique discussed in \autoref{hyperreduction} is set to $\tau=0.01$.
The ability of the proposed \emph{HP-VpROM} to accurately infer response fields that are relevant for dynamic structural systems, while providing accelerated model evaluations is exhibited herein.

\begin{figure*}[!ht]
        \begin{subfigure}[b]{0.49\textwidth}   
            \centering 
            \includegraphics[scale=0.60]{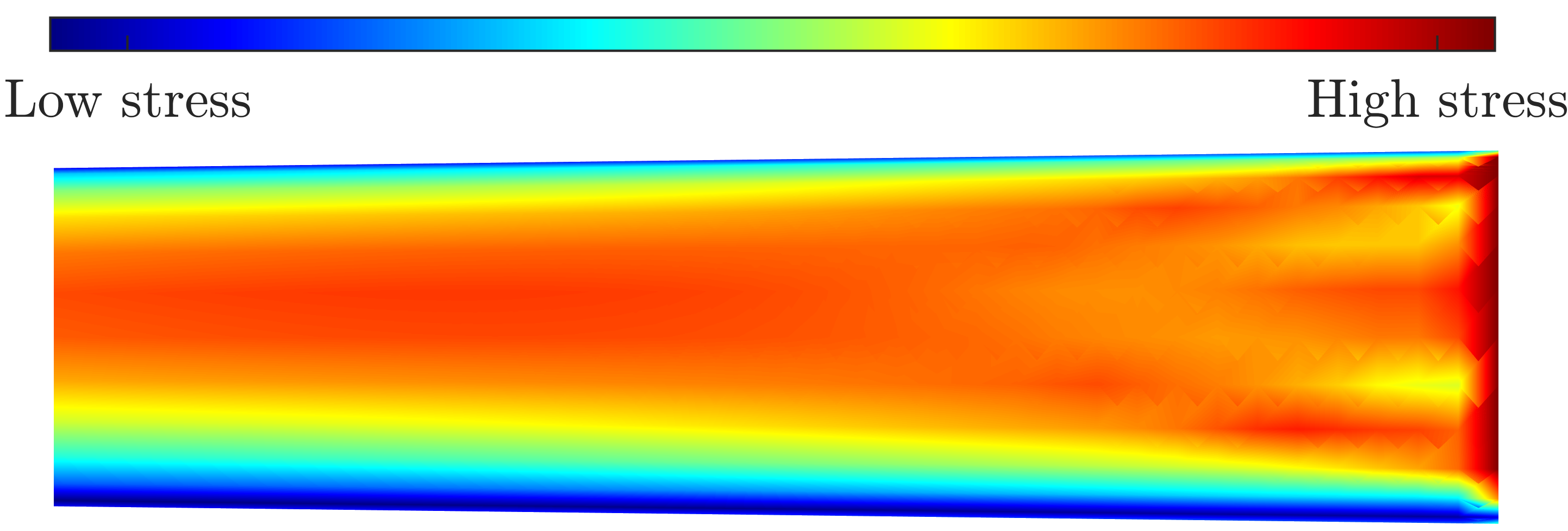}
            \caption[]%
            {\ac{FOM} stress visualisation for $[3.40, 423.65, 1.01]$.}    
            \label{fig:stressHFMa}
        \end{subfigure}
        \quad
        \begin{subfigure}[b]{0.49\textwidth}   
            \centering 
            \includegraphics[scale=0.60]{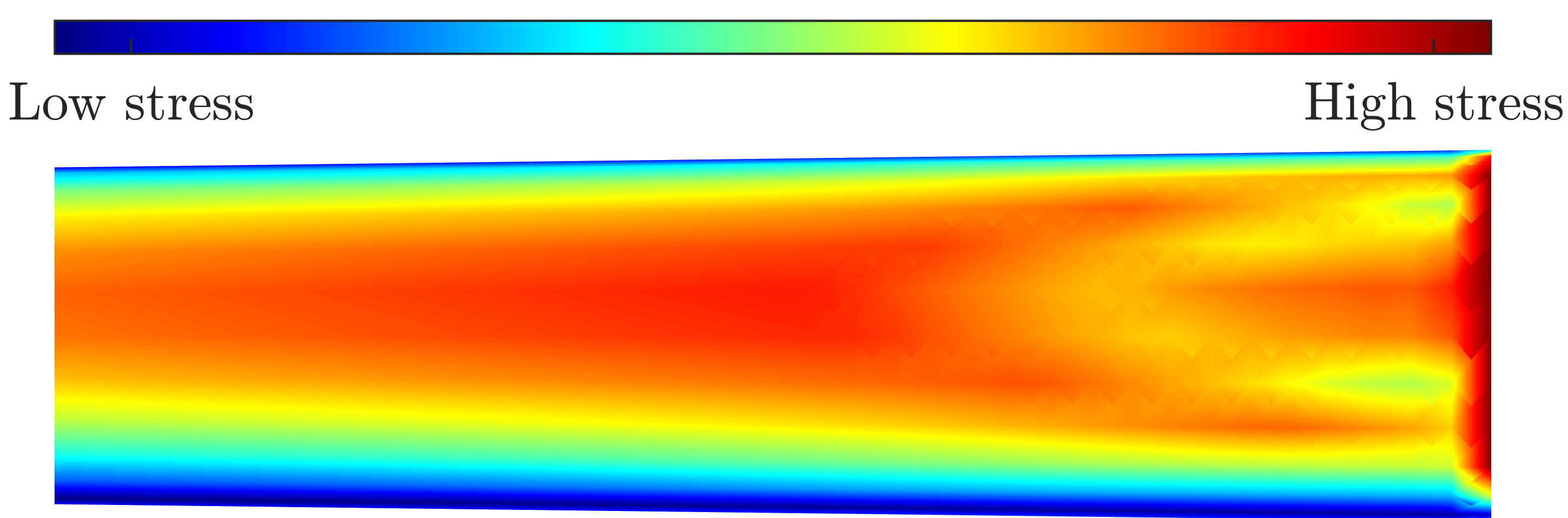}
            \caption[]%
            {\emph{VpROM} stress visualisation for $[3.40, 423.65, 1.01]$.} 
            \label{fig:stressROMa}
        \end{subfigure}
         \begin{subfigure}[b]{0.49\textwidth}   
            \centering 
            \includegraphics[scale=0.60]{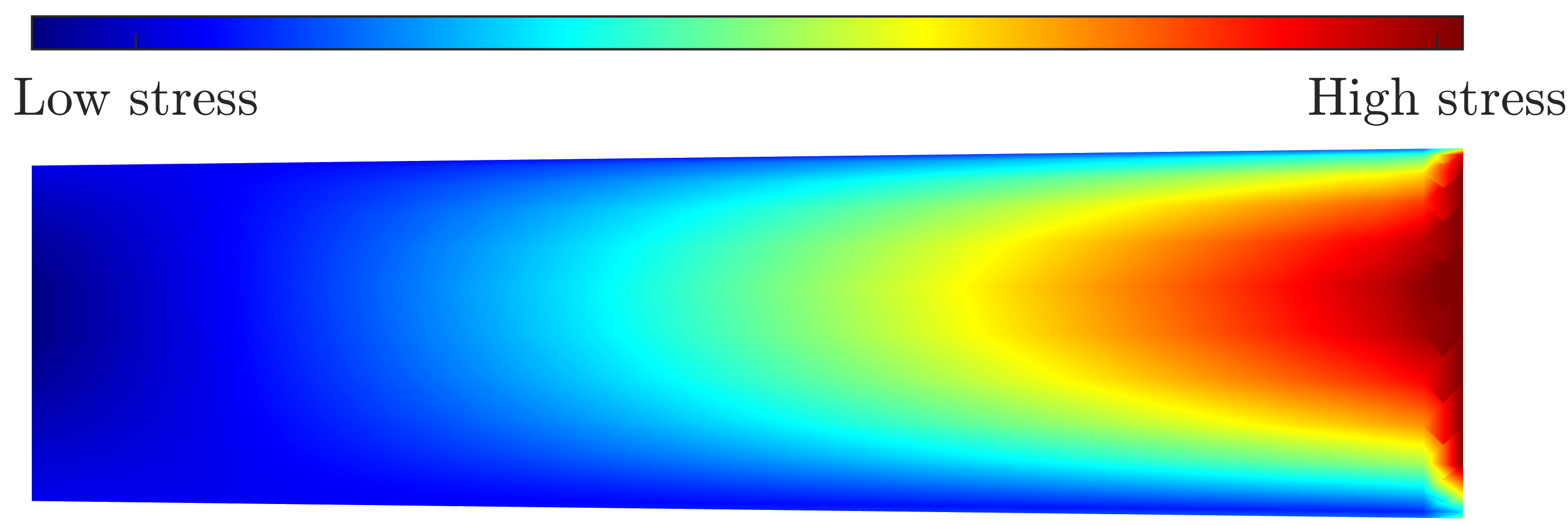}
            \caption[]%
            {\ac{FOM} stress visualisation for $[3.49, 411.40, 1.17]$.}    
            \label{fig:stressHFMb}
        \end{subfigure}
        \quad
        \begin{subfigure}[b]{0.49\textwidth}   
            \centering 
            \includegraphics[scale=0.60]{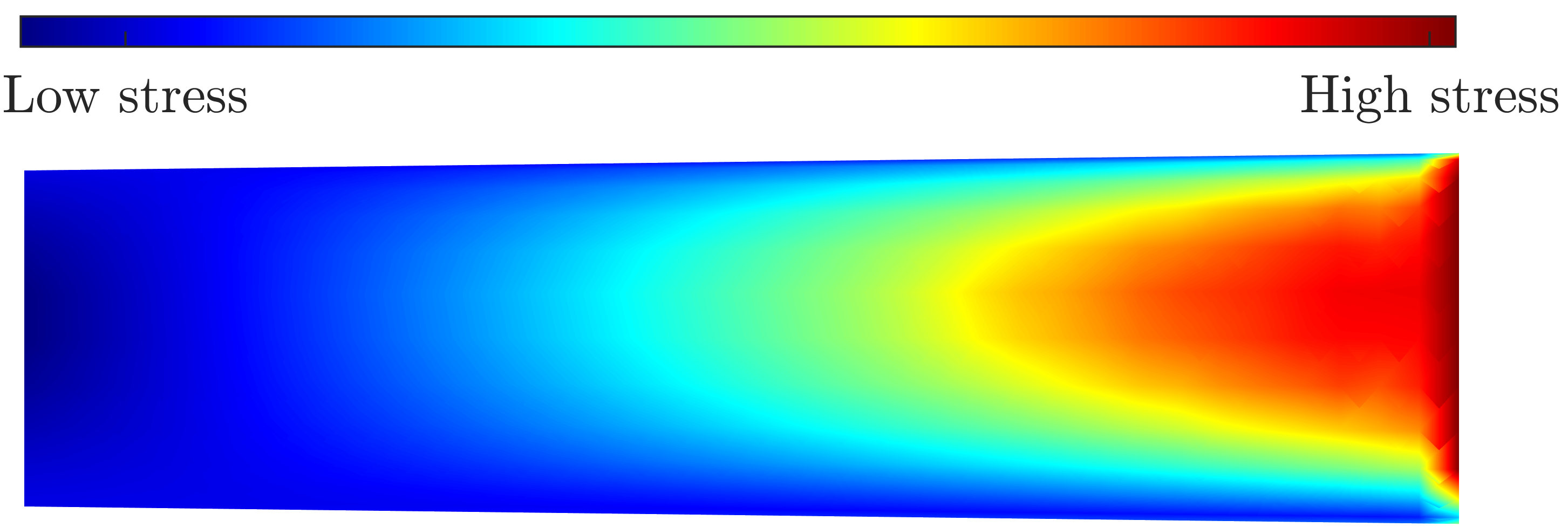}
            \caption[]%
            {\emph{VpROM} stress visualisation for $[3.49, 411.40, 1.17]$} 
            \label{fig:stressROMb}
        \end{subfigure}
        \caption{Visualisation of the approximation of internal stresses achieved by the proposed \emph{VpROM} using nodal averaging. Only the yielding domain is visualised, which extends to one-third of the total height. } 
        \label{fig:stresses}
\end{figure*}

In \autoref{fig:stresses} a visualisation of the \emph{HP-VpROM} approximation for the internal stress field in the yielding domain is provided for two validation examples, which feature different dynamic behavior.
The respective high-fidelity field is also visualised via the \ac{FOM} for reference purposes. Stresses and strains are important metrics to be monitored in many structural applications; thus, the ability to capture their distribution accurately is often of critical importance.
Despite the minor discrepancies observed, the overall quality of the \emph{HP-VpROM} approximation illustrated in \autoref{fig:stresses} indicates an effective low-order representation, able to deliver high precision estimates of stress state distributions.
This exemplifies the potential utility of the proposed \emph{HP-VpROM} in condition monitoring, fatigue, or damage localisation.

A more comprehensive summary of all aspects of the \ac{ROM}s' performance is provided in \autoref{tab:ErrorsMonopile}. 
Specifically, the average and maximum error measures of the respective approximations on displacement and acceleration time histories are summarised for all hyper-reduced surrogates of \autoref{tab:pROMS}.
In addition, the precision is reported for two example validation samples in the extreme regions of the input domain for reference.

Similar to what was already observed in the previous case study in \autoref{example:BWlinks}, the hyper-reduced variant of the \emph{VpROM}, namely \emph{HP-VpROM}, delivers a superior surrogate in terms of accuracy.
The respective low maximum $err_u$ measure on capturing the displacements indicates a robust approximation, whereas the other two \ac{ROM}s experience performance outliers where accuracy deteriorates significantly.
At the same time, \emph{HP-VpROM} achieves an average discrepancy lower than $1\%$, implying a high precision representation.
Similar conclusions can be drawn by observing the respective measures on the two validation samples offered as additional examples. 
Regarding the inference of the acceleration response, the \emph{VpROM} maintains its superior accuracy, although it seems only marginally better than the two alternative \ac{ROM}s implemented.

\begin{table}[!htb]
	\caption{Performance measures for the hyper-reduced \ac{ROM}s from \autoref{tab:pROMS}. Two validation samples are presented as an example, along with the respective median and max error measures from \autoref{eq:errors}. Efficiency is also reported.}
	\label{tab:ErrorsMonopile}
		\begin{tabular} { p{2.5cm} | p{0.8cm} p{0.8cm} p{0.8cm} p{0.8cm} | p{0.8cm} p{0.8cm} p{0.8cm} p{0.8cm} | p{1.5cm} p{1.2cm} }
        \hline
		  & \multicolumn{2}{c}{Sample}
        & \multicolumn{2}{c}{Sample}
        & \multicolumn{2}{c}{Average}
        & \multicolumn{2}{c}{Maximum}
        & \multicolumn{2}{c}{Efficiency}
        \\
        & \multicolumn{2}{c}{[2.5,375,0.80]}
        & \multicolumn{2}{c}{[3.75,450,1.2]}
        & \multicolumn{2}{c}{Error}
        & \multicolumn{2}{c}{Error}
        & Speed-Up
        & CPU
        \\
        &$err_u$&$err_{\ddot{u}}$&$err_u$ & $err_{\ddot{u}}$ & $err_u$ & $err_{\ddot{u}}$ & $err_u$ & $err_{\ddot{u}}$ & Factor & timing
        \\
		\hline
        \ac{FOM} & - & - & - & - & - & - & - & - & 1.00 & 3990 (s) \\
		HP-\ac{MACpROM} & 2.01\%& 6.44\%& 1.54\% & 5.88\% & 2.34$\%$ & 6.24$\%$ & 8.89$\%$ & 7.64$\%$ & 35.40 & 112 (s) \\
		HP-CpROM & 1.44\% & 6.48\% & 1.06\% & 5.87\% & 2.05$\%$ & 6.22$\%$ & 5.09$\%$ & 7.48$\%$ & 36.69 & 109 (s) \\
		HP-VpROM & 0.51\%& 6.44\%& 0.84\%& 5.89\% & 0.91$\%$ & 5.98$\%$ & 1.88$\%$ & 7.13$\%$ & 35.67 & 112 (s) \\
        \hline
	\end{tabular}
\end{table}

The utility of the proposed \emph{HP-VpROM} for applications, in which (near) real-time model evaluations are required is also documented in \autoref{tab:ErrorsMonopile}.
The hyper-reduction technique, along with the ability of the \ac{ROM} to propagate the dynamics in a proper low-order subspace achieves a substantial computational toll reduction and accelerated computations. 
The respective average speed-up factor $t_{FOM}/t_{ROM}$ reported in \autoref{tab:ErrorsMonopile} implies significant savings in computational resources during model evaluations.

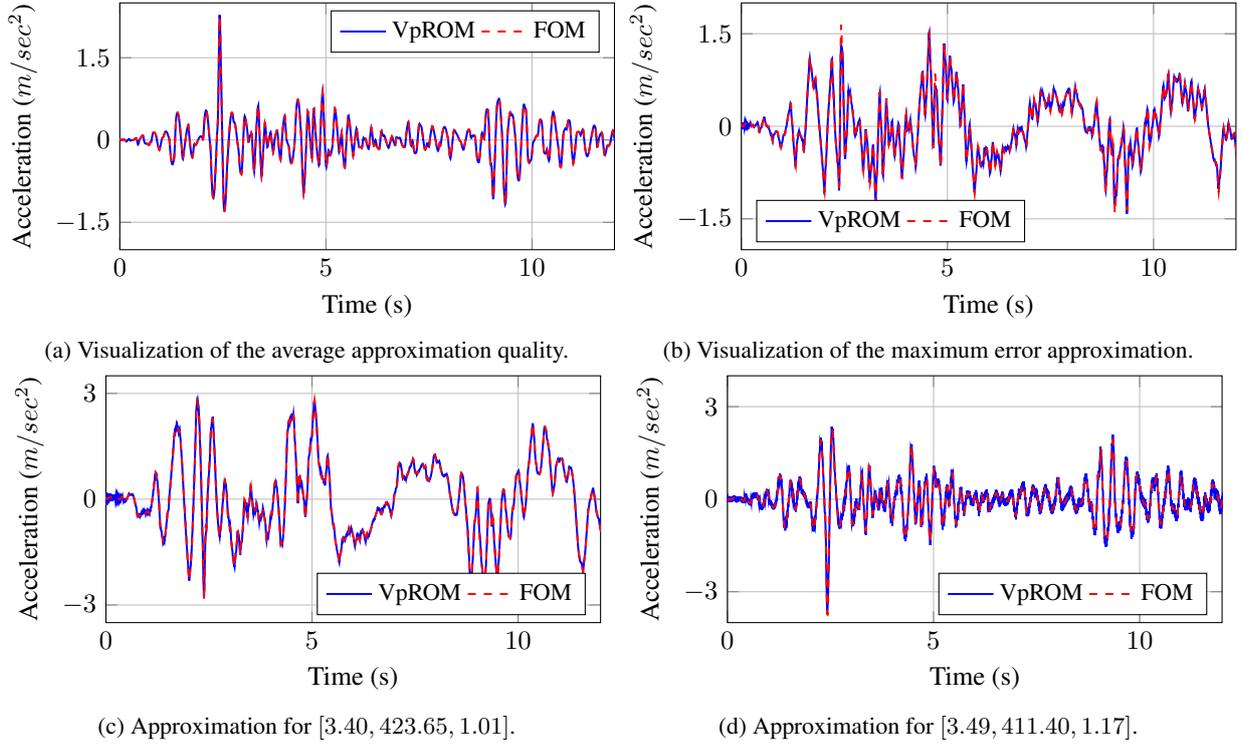
\begin{figure}[!htb]
\centering
    \begin{subfigure}[c]{0.49\textwidth}
        \centering
        \begin{tikzpicture}
	\begin{axis}[
		name = response,
		xmin = 0,
		xmax = 12.0,
		ymin = -2.0,
		ymax = 2.5,
		xtick = {0, 5,10,15},
		ytick = {-1.5,0,1.5},
		xlabel = {Time (s)},
		ylabel = {Acceleration ($m/sec^2$)},
		grid = both,
		width=\textwidth,
		height=0.60\textwidth,
        legend style={legend columns=-1, legend pos=north east,nodes={font=\fontsize{\figureFontSize pt}{\figureFontSize pt}\selectfont}}
		]
			
		\addplot[color = blue , line width=0.75pt] table[x=time, y=two] {figures_tikz/ResponseOverview.dat};

        \addplot[color = red, dashed, line width=0.75pt] table[x=time, y=one] {figures_tikz/ResponseOverview.dat};
  		
		\addlegendentry{VpROM}
        \addlegendentry{FOM}
	\end{axis}
	
\end{tikzpicture}
        \captionsetup{skip=-6pt}
        \caption{Visualization of the average approximation quality.} 
    \end{subfigure}     
    \begin{subfigure}[c]{0.49\textwidth}
        \centering
        \begin{tikzpicture}
	\begin{axis}[
		name = response,
		xmin = 0,
		xmax = 12.0,
		ymin = -2.0,
		ymax = 2.0,
		xtick = {0, 5,10,15},
		ytick = {-1.5,0,1.5},
		xlabel = {Time (s)},
		ylabel = {Acceleration ($m/sec^2$)},
		grid = both,
		width=\textwidth,
		height=0.60\textwidth,
        legend style={legend columns=-1, legend pos=south west,nodes={font=\fontsize{\figureFontSize pt}{\figureFontSize pt}\selectfont}}
		]
			
		\addplot[color = blue , line width=0.75pt] table[x=time, y=four] {figures_tikz/ResponseOverview.dat};

        \addplot[color = red, dashed, line width=0.75pt] table[x=time, y=three] {figures_tikz/ResponseOverview.dat};
  		
		\addlegendentry{VpROM}
        \addlegendentry{FOM}
	\end{axis}
	
\end{tikzpicture}
        \captionsetup{skip=-6pt}
        \caption{Visualization of the maximum error approximation.}
    \end{subfigure}
     \begin{subfigure}[c]{0.49\textwidth}
        \centering
        \begin{tikzpicture}
	\begin{axis}[
		name = response,
		xmin = 0,
		xmax = 12.0,
		ymin = -3.5,
		ymax = 3.5,
		xtick = {0, 5,10,15},
		ytick = {-3.0,0,3.0},
		xlabel = {Time (s)},
		ylabel = {Acceleration ($m/sec^2$)},
		grid = both,
		width=\textwidth,
		height=0.60\textwidth,
        legend style={legend columns=-1, legend pos=south east,nodes={font=\fontsize{\figureFontSize pt}{\figureFontSize pt}\selectfont}}
		]
			
		\addplot[color = blue , line width=0.75pt] table[x=time, y=six] {figures_tikz/ResponseOverview.dat};

        \addplot[color = red, dashed, line width=0.75pt] table[x=time, y=five] {figures_tikz/ResponseOverview.dat};
  		
		\addlegendentry{VpROM}
        \addlegendentry{FOM}
	\end{axis}
	
\end{tikzpicture}
        \caption{Approximation for $[3.40, 423.65, 1.01]$.} 
    \end{subfigure}     
    \begin{subfigure}[c]{0.49\textwidth}
        \centering
        \begin{tikzpicture}
	\begin{axis}[
		name = response,
		xmin = 0,
		xmax = 12.0,
		ymin = -4.0,
		ymax = 4.0,
		xtick = {0, 5,10,15},
		ytick = {-3.0,0,3.0},
		xlabel = {Time (s)},
		ylabel = {Acceleration ($m/sec^2$)},
		grid = both,
		width=\textwidth,
		height=0.60\textwidth,
        legend style={legend columns=-1, legend pos=south east,nodes={font=\fontsize{\figureFontSize pt}{\figureFontSize pt}\selectfont}}
		]
			
		\addplot[color = blue , line width=0.75pt] table[x=time, y=eight] {figures_tikz/ResponseOverview.dat};

        \addplot[color = red, dashed, line width=0.75pt] table[x=time, y=seven] {figures_tikz/ResponseOverview.dat};
  		
		\addlegendentry{VpROM}
        \addlegendentry{FOM}
		
	\end{axis}
	
\end{tikzpicture}
        \caption{Approximation for $[3.49, 411.40, 1.17]$.}
    \end{subfigure}
	\caption{Visualisation of the different levels of the approximation quality achieved using the \emph{HP-VpROM}. The estimation is reported for various response patterns the system exhibits depending on its parametric features.}
	\label{fig:RespOverviewMonopile}
\end{figure}

The reported performance measures highlight the fact that the proposed \ac{ROM} framework remains robust and precise, also when coupled with hyper-reduction.
The generative model injected guarantees the ability of the \emph{HP-VpROM} to capture different dynamic trends in the response and avoid accuracy outliers.
This is indicatively visualised in \autoref{fig:RespOverviewMonopile}, where the acceleration response of the full-order model and the respective \emph{HP-VpROM} approximation is illustrated for validation samples located near the edges of the input domain.

The different patterns in the system's behavior are demonstrated clearly, along with the ability of the assembled \emph{HP-VpROM} to capture the different trends sufficiently accurately.
Despite the minor discrepancies observed, especially when high-frequency components are present in the response like in the bottom right example, the \emph{HP-VpROM} maintains a robust performance with a high-quality approximation across the input domain.

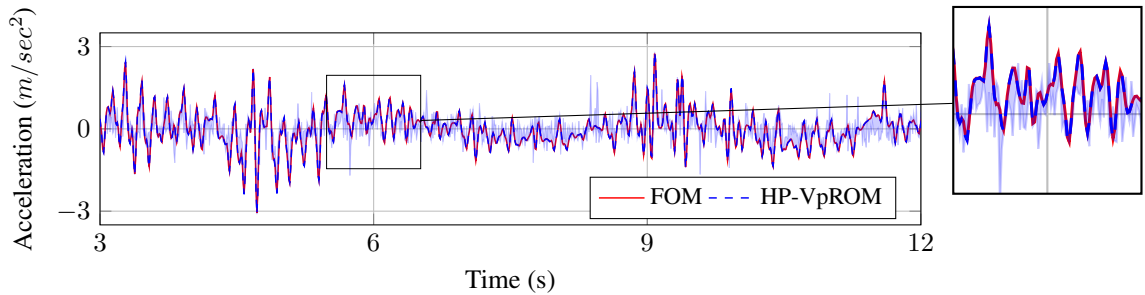
\begin{figure}[!htb]
	\centering
\begin{tikzpicture}[spy using outlines=
	{rectangle, magnification=2, anchor = center, connect spies}]
	\begin{axis}[
		name = response,
		xmin = 3,
		xmax = 12,
		ymin = -3.5,
		ymax = 3.5,
		xtick = {3,6,9,12,15},
		ytick = {-3.0,0,3.0},
		xlabel = {Time (s)},
		ylabel = {Acceleration ($m/sec^2$)},
		grid = both,
		width=0.75\textwidth,
		height=0.25\textwidth,
		legend style={legend columns=-1, legend pos=south east,nodes={font=\fontsize{\figureFontSize pt}{\figureFontSize pt}\selectfont}}
		]
		
		\addplot[name path=F, color = red, line width=0.55pt] table[x=time, y=FOM] {figures_tikz/UQplot_Acc.dat};
		\addplot[color = blue, dashed, line width=0.55pt] table[x=time, y=ROM] {figures_tikz/UQplot_Acc.dat};
		\addplot[name path=G, color = blue,opacity=0.2] table[x=time, y=Max] {figures_tikz/UQplot_Acc.dat};			
		\tikzfillbetween[of=F and G]{blue, opacity=0.2};
		\addplot[name path=K, color = blue,opacity=0.2] table[x=time, y=Inf] {figures_tikz/UQplot_Acc.dat};			
		\tikzfillbetween[of=F and K]{blue, opacity=0.2};
		
		\addlegendentry{FOM}
		\addlegendentry{HP-VpROM}

		\coordinate (spypoint) at (axis cs:6,0.25);
		
	\end{axis}
	
	\path node[anchor= center] (magnifyglass) at (0.755\textwidth,0.10\textwidth) {};
	
	\spy [black, size=0.15\textwidth] on (spypoint)
	in node[fill=white] at (magnifyglass);
\end{tikzpicture}
\caption{Quality and confidence bounds of the \emph{HR-VpROM} approximation. Evaluation is performed on the degree of freedom with the maximum absolute response. The shaded area quantifies the uncertainty of the response inference.}
 \label{fig:THsmon}
\end{figure}

The coupling of the \ac{ROM} with a \ac{cVAE}-based generative model offers an additional feature, namely the quantification of the uncertainty in the estimations.
Relying on the probabilistic nature of the latent space of the assembled \ac{cVAE} the proposed \emph{HP-VpROM} comes with confidence bounds on its predictions.
This is exemplified in \autoref{fig:THsmon}, where the approximation for the acceleration is visualised for a representative validation sample.
The shaded region in \autoref{fig:THsmon} represents the uncertainty of the respective prediction and may be used as a confidence measure during when using the \ac{ROM} predictions this improves the utility of the scheme compared to determinstic methods.

\section{Limitations and concluding remarks}
\label{Conclusions}
This work demonstrates the use of a \ac{cVAE}-based \ac{ROM}, termed as \emph{VpROM}, as an extension to state of the art methods for generating local reduction bases for nonlinear parametric ROMs. 

The following conclusions are drawn:
\begin{itemize}
    \item cVAE neural networks can successfully be used to generate local bases for nonlinear parametric ROMs with high dimensional parameterisation and strongly nonlinear behaviour. 
    \item The verification of the proposed scheme on a large-scale system results in significantly accelerated model evaluations, almost 40 times faster compared to the \ac{FOM}.
    \item The cVAE can outperform current state of the art methods such as interpolation and clustering algorithms in terms of precision of the R\ac{ROM}.
    \item The \emph{VpROM} formulation offers the additional benefit of encoding the uncertainty in the predicted local bases and the ability to propagate this in the predicted response.

\end{itemize}

The newly developed method is demonstrated on two simulated nonlinear systems, which are parameterised in terms of both system and loading traits. The first example demonstrates the viability of the method for a system of high dimensional parametric dependency exhibiting strongly varying nonlinear behaviour. The second example verifies the proposed scheme on a large-scale system, with the inclusion of hyper-reduction, in order to demonstrate the utility of the method for hyper-accelerated model evaluations. In both examples, the potential of the method for quantifying uncertainty on its estimates is also demonstrated.

The main limitation of the method in comparison to current methodologies is the relative complexity of the training process of the \ac{cVAE} models. Owing to their very flexible nature neural network methods, such as VAEs, comprise a relatively high number of hyper-parameters that must be tuned for optimal performance. When training the \emph{VpROM}, it is necessary to select such hyper-parameters, which include  for instance, the number of layers in the network and the number of neurons in each of these layers, the activation functions used in the network, and the learning rate of the optimisation algorithm. There exist a number of heuristics for choosing such parameters, such as grid search, yet it remains worth highlighting that this process requires more effort than is required for other state-of-the-art methods exploiting clustering or interpolation methods.

\section*{Data Availability }
The data that support the findings of this study are available from the authors, \textbf{TS}, \textbf{KV}, upon reasonable request.

\section*{Declaration of Competing Interest}
The authors declare that they have no known competing financial interests or personal relationships that could have appeared to influence the work reported in this paper.

\section*{Acknowledgements}

This research has received funding by the Sandia National Laboratories. Sandia National Laboratories is a multimission laboratory managed and operated by National Technology and Engineering Solutions of Sandia, LLC, a wholly owned subsidiary of Honeywell International Inc., for the U.S. Department of Energy’s National Nuclear Security Administration under contract DE-NA0003525.  This paper describes objective technical results and analysis. Any subjective views or opinions that might be expressed in the paper do not necessarily represent the views of the U.S. Department of Energy or the United States Government. Further funding was received from the European Union’s Horizon 2020 research and innovation programme under the Marie Skłodowska-Curie grant agreement No 764547 and from the EPSRC under grant agreement EP/R004900/1.

\bibliographystyle{plain}


\section*{CRediT authorship contribution statement}
\textbf{Thomas Simpson}: Conceptualization, Methodology, Software, Writing. \textbf{Konstantinos Vlachas}: Conceptualization, Methodology, Validation, Writing. \textbf{Anthony Garland}: Resources, Writing - review \& editing. \textbf{Nikolaos Dervilis}: Writing - review \& editing, Supervision. \textbf{Eleni Chatzi}: Conceptualization, Resources, Writing - review \& editing, Supervision.

\end{document}